\documentclass[cleveref,a4paper,UKenglish]{lipics-v2021}
\usepackage{cite} 
\usepackage{amsmath}
\usepackage{amssymb}  
\usepackage{xspace}

\usepackage{tikz}
\usepackage{todonotes}
\tikzstyle{every picture} = [>=latex]
\usetikzlibrary{calc,arrows,decorations.markings}

\def\ca#1{{\cal#1}}

\hideLIPIcs
\nolinenumbers

\title{Twin-width of Planar Graphs is at most~9, and at most~6 when Bipartite Planar}

\author{Petr Hlin{\v e}n\'y}{Masaryk University, Brno,
  Czech republic}{hlineny@fi.muni.cz}{https://orcid.org/0000-0003-2125-1514}{}
  
\authorrunning{P.\ Hlin\v{e}n\'y}
\Copyright{Petr Hlin\v{e}n\'y}


\funding{Supported by the {Czech Science Foundation}, project no.~{20-04567S}.}


\begin{document}
\maketitle

\begin{abstract}
The structural parameter twin-width was introduced by Bonnet et al.~in~[FOCS 2020], 
and already this first paper included an asymptotic argument bounding the twin-width of planar graphs by a non-explicit constant.
Quite recently, we have seen first small explicit upper bounds of $183$ by Jacob and Pilipczuk [arXiv, January 2022], $583$ by Bonnet et al.~[arXiv, February 2022],
and of $37$ by Bekos et al. [arXiv, April 2022].
We prove that the twin-width of planar graphs is at most~$9$.
Furthermore, if a planar graph is also bipartite, then its twin-width is at most~$6$.
\end{abstract}

\section{Introduction}

Twin-width is a new structural width measure of graphs introduced in 2020 by Bonnet, Kim, Thomass\'e and Watrigant~\cite{DBLP:journals/jacm/BonnetKTW22}.

Consider only simple graphs in the coming definition.\footnote{In general, the concept of twin-width is defined for binary relational structures of a finite signature, and
so one may either define the twin-width of a multigraph as the twin-width of its simplification, or allow only bounded multiplicities of edges and use the more general matrix definition of twin-width.}
A \emph{trigraph} is a simple graph $G$ in which some edges are marked as {\em red}, and we then naturally speak about \emph{red neighbours} and \emph{red degree} in~$G$. 
For a pair of (possibly not adjacent) vertices $x_1,x_2\in V(G)$, we define a \emph{contraction} of the pair $x_1,x_2$ as the operation
creating a trigraph $G'$ which is the same as $G$ except that $x_1,x_2$ are replaced with a new vertex $x_0$ (said to {\em stem from $x_1,x_2$}) such that:
\begin{itemize}
\item 
the full neighbourhood of $x_0$ in $G'$ (i.e., including the red neighbours), denoted by $N_{G'}(x_0)$, 
equals the union of the neighbourhoods of $x_1$ and $x_2$ in $G$ except $x_1,x_2$ themselves, that is,
$N_{G'}(x_0)=(N_G(x_1)\cup N_G(x_2))\setminus\{x_1,x_2\}$, and
\item 
the red neighbours of $x_0$, denoted by $N_{G'}^r(x_0)$, inherit all red neighbours of $x_1$ and of $x_2$ and add those in $N_G(x_1)\Delta N_G(x_2)$,
that is, $N_{G'}^r(x_0)=\big(N_{G}^r(x_1)\cup N_G^r(x_2)\cup(N_G(x_1)\Delta N_G(x_2))\big)\setminus\{x_1,x_2\}$, where $\Delta$ denotes the~symmetric~difference.
\end{itemize}
A \emph{contraction sequence} of a trigraph $G$ is a sequence of successive contractions turning~$G$  into a single vertex, and its \emph{width} is the maximum red degree of any vertex in any trigraph of the sequence.
The \emph{twin-width} of a trigraph $G$ is the minimum width over all possible contraction sequences of~$G$.

\smallskip
As shown already in the pioneering paper on this concept~\cite{DBLP:journals/jacm/BonnetKTW22}, the twin-width of planar graphs is bounded
(but there was no explicit number given there).
In the first half of 2022, first explicit upper bounds on the twin-width of planar graphs have appeared;
$183$ by Jacob and Pilipczuk~\cite{DBLP:journals/corr/abs-2201-09749},
$583$ by Bonnet et al.~\cite{DBLP:journals/corr/abs-2202-11858} (this paper more generally bounds the twin-width of so-called $k$-planar graphs with an asymptotic exponential function of $O(k)$),
and most recently $37$ by Bekos et~al.~\cite{DBLP:journals/corr/abs-2204-11495} (which more generally bounds the twin-width of so-called $h$-framed graphs,
and in particular gives also an upper bound of $80$ on the twin-width of $1$-planar graphs).

It is worth to mention that all three papers~\cite{DBLP:journals/corr/abs-2201-09749,DBLP:journals/corr/abs-2202-11858,DBLP:journals/corr/abs-2204-11495}, 
more or less explicitly, use the product structure machinery of planar graphs (cf.~\cite{DBLP:journals/jacm/DujmovicJMMUW20}).
We give the following improved bounds on the twin-width of planar graphs:

\begin{theorem}\label{thm:twwplanar}
The twin-width of any simple planar graph is at most~$9$, and a corresponding contraction sequence can be found in linear time.
\end{theorem}
\begin{theorem}\label{thm:twwbiplanar}
The twin-width of any simple {\em bipartite} planar graph is at most~$6$, and a corresponding contraction sequence can be found in linear time.
\end{theorem}
Our proofs do not use the product structure machinery, but another recursive decomposition of planar graphs which is fine-tuned for the purpose of obtaining the results.

\section{Proving the bound for general planar graphs}

We start with a few technical definitions needed for the proof of \Cref{thm:twwplanar}.

Let $G$ be a graph and $r\in V(G)$ a fixed vertex.
The \emph{BFS layering of $G$ determined by~$r$} is the vertex partition $\ca L=(L_0=\{r\},L_1,L_2,\ldots)$ of $G$ such that 
$L_i$ contains all vertices of $G$ at distance exactly $i$ from $r$.
A path $P\subseteq G$ is {\em $r$-geodesic} if $P$ is a subpath of some shortest path from $r$ to any vertex of~$G$
(in particular, $P$ intersects every layer of $\ca L$ in at most one vertex).
Let $T$ be a {\em BFS tree} of $G$ rooted at the vertex $r\in V(G)$ as above
(that is, for every vertex $v\in V(G)$, the distance from $v$ to $r$ is the same in $G$ as in $T$).
A path $P\subseteq G$ is {\em vertical (wrt.~implicit~$T$)} if $P$ is a subpath of some root-to-leaf path of~$T$.
Notice that a vertical path is $r$-geodesic, but the converse may not be true.

Observe the following:
\begin{claim}\label{cl:layeri3}
For every edge $\{v,w\}$ of $G$ with $v\in L_i$ and $w\in L_j$, we have $|i-j|\leq 1$, and so
a contraction of a pair of vertices from $L_i$ may create new red edges only to the remaining vertices of~$L_{i-1}\cup L_i\cup L_{i+1}$.
\end{claim}

Consider a BFS layering $\ca L$ of $G$ as above.
A \emph{partial contraction sequence} of $G$ is defined in the same way as a contraction sequence of $G$, except that it does not have to end with a single-vertex graph.
Such a sequence is \emph{$\ca L$-respecting} if every step contracts only pairs belonging to $L_j\in\ca L$ for some~$j$.
For reference, when $G'$ is a trigraph resulting from an $\ca L$-respecting partial contraction sequence of $G$, 
we will write $L_j[G']$ for the vertex set that stems from $L_j$ by these contractions (hence, $L_j=L_j[G]$). 

We will deal with {\em plane graphs}, which are planar graphs with a given (combinatorial) embedding in the plane, and one marked {\em outer face}
(the remaining faces are then {\em bounded}).
A plane graph is a {\em plane triangulation} if every face of its embedding is a triangle.
Likewise, a plane graph is a {\em plane quadrangulation} if every face of its embedding is of length~$4$.

In this planar setting, consider now a plane graph $G$, and a BFS tree $T$ spanning $G$ and rooted in a vertex $r$ of the outer face of $G$,
and picture (for clarity) the embedding $G$ such that $r$ is the vertex of $G$ most at the top.
For two adjacent vertices $u,v\in V(G)$, $uv\in E(G)$, we say that {\em$u$ is to the left of~$v$} (wrt.~$T$) if neither of $u,v$
lies on the vertical path from $r$ to the other, and the following holds;
if $r'$ is the least common ancestor of $u$ and $v$ in $T$ and $P_{r',u}$ (resp., $P_{r',v}$) denote the vertical path from $r'$ to $u$ (resp., $v$),
then the cycle $(P_{r',u}\cup P_{r',v})+uv$ has the triple $(r',u,v)$ in this counter-clockwise cyclic order.

A BFS tree $T$ of $G$ with the BFS layering $\ca L=(L_0,L_1,\ldots)$ is called {\em left-aligned} if there is no edge $f=uv$ of $G$ such that, for some index $i$,
$u\in L_{i-1}$ and $v\in L_i$, and $u$ is to the left of~$v$
(an informal meaning is that one cannot choose another BFS tree of $G$ which is ``more to the left'' of~$T$ in the geometric picture of $G$ and~$T$, such as by picking the edge $uv$ instead of the parental edge of $v$ in~$T$).

\begin{claim}\label{clm:existslal}
Given a simple plane graph $G$, and a vertex $r$ on the outer face, there exists a left-aligned BFS tree of $G$
and can be found in linear time.
\end{claim}

\begin{claimproof}
For this proof, we have to extend the above relation of ``being left of'' to edges emanating from a common vertex of $G$.
So, for an arbitrary BFS tree T of $G$ and edges $f_1,f_2\in E(G)$ incident to $v\in V(G)$, such that neither of $f_1,f_2$ is the parental edge of $v$ in~$T$,
we write $f_1\leq_lf_2$ if there exist adjacent vertices $u_1,u_2\in V(G)$ such that $u_1$ is to the left of~$u_2$,
the least common ancestor of $u_1$ and $u_2$ in $T$ is~$v$ and, for $i=1,2$, the edge $f_i$ lies on the vertical path from $u_i$ to~$v$.
Observe the following; if $f_0$ is the parental edge of $v$ in~$T$ (or, in case of $v=r$, $f_0$ is a ``dummy edge'' pointing straight up from~$r$),
then $f_1\leq_lf_2$ implies that the counter-clockwise cyclic order around $v$ is $(f_0,f_1,f_2)$.
In particular, $\leq_l$ can be extended into a linear order on its domain.

We first run a basic linear-time BFS search from $r$ to determine the BFS layering $\ca L$~of~$G$.
Then we start the construction of a left-aligned BFS tree $T\subseteq G$ from $T:=\{r\}$, and we recursively (now in a ``DFS manner'') proceed as follows:
\begin{itemize}
\item Having reached a vertex $v\in V(T)\subseteq V(G)$ such that $v\in L_i$, and denoting by $X:=(N_G(v)\cap L_{i+1})\setminus V(T)$
all neighbours of $v$ in $L_{i+1}$ which are not in $T$ yet, we add to $T$ the nodes $X$ and the edges from $v$ to~$X$.
\item We order the vertices in $X$ using the cyclic order of edges emanating from $v$ to have it compatible with $\leq_l$ at~$v$,
and in this increasing order we recursively (depth-first, to be precise) call the procedure for them.
\end{itemize}
The result $T$ is clearly a BFS tree of $G$. Assume, for a contradiction, 
that $T$ is not left-aligned, and let $u_1\in L_{i-1}$ and $u_2\in L_i$ be a witness pair of it, where $u_1u_2\in E(G)$ and $u_1$ is to the left of~$u_2$.
Let $v$ be the least common ancestor of $u_1$ and $u_2$ in $T$, and let $v_1$ and $v_2$ be the children of $v$ on the $T$-paths from $v$ to $u_1$ and $u_2$, respectively.
So, by the definition, $vv_1\leq_lvv_2$ at $v$, and hence when $v$ has been reached in the construction of $T$,
its child $v_1$ has been taken for processing before the child $v_2$.
Consequently, possibly deeper in the recursion, $u_1$ has been processed before the parent of~$u_2$ and,
in particular, the procedure has added the edge $u_1u_2$ into $T$, a contradiction to $u_1$ being to the left of~$u_2$.

This recursive computation is finished in linear time, since every vertex of $G$ is processed only in one branch of the recursion,
and one recursive call takes time linear in the number of incident edges (to~$v$).
\end{claimproof}

For a plane graph $G$ and its cycle $C$, the {\em subgraph of $G$ bounded by $C$}, denoted by $G_C$, 
is the subgraph of $G$ formed by the vertices and edges of $C$ and the vertices and edges of $G$ drawn inside~$C$
-- formally, in the region of the plane bounded by $C$ and not containing the outer face.
It is easy to turn any planar graph into a simple plane triangulation by adding vertices and incident edes into each non-triangular face.
Furthermore, twin-width is non-increasing when taking induced subgraphs, and so it suffices to focus on plane triangulations
when stating the following core lemma:
\begin{lemma}\label{lem:core}
Let $G$ be a plane triangulation and $T$ be a left-aligned BFS tree of $G$ rooted at a vertex $r\in V(G)$ of the outer triangular face 
and defining the BFS layering $\ca L=(L_0=\{r\},L_1,L_2,\ldots)$ of~$G$.
Assume that $P_i$, $i=1,2$, are two edge-disjoint vertical paths of $G$ of lengths at least $1$,
where $P_i$ has the ends $u_i$ and $v_i$ such that $u_1=u_2\in L_\ell$ is the vertex of $P_1\cup P_2$ closest to~$r$ (i.e., at distance $\ell$ from~$r$),
and $f=\{v_1,v_2\}\in E(G)$ is an edge such that $v_1$ is to the left of $v_2$.
Then $C:=(P_1\cup P_2)+f$ is a cycle of $G$, and let $G_C$ be the subgraph of $G$ bounded by~$C$,
let $U:=V(G_C)\setminus V(C)$ and $W:=V(G)\setminus V(G_C)$ (so, $V(G)$ is tri-partitioned into $U$, $W$ and~$V(C)$,
and $U\cap L_j=\emptyset$ for all $j\leq\ell$).

Then there exists an $\ca L$-respecting partial contraction sequence of $G$ which contracts only pairs of vertices that are in or stem from $U$,
results in a trigraph $G^0$, and satisfies the following conditions for every trigraph $G'$ along this contraction sequence from $G$ to~$G^0$:
\begin{enumerate}[(I)]
\item For $U':=V(G')\setminus(V(C)\cup W)$ (which are the vertices that are in or stem from $U$ in~$G'$),
every vertex of $U'$ in $G'$ has red degree at most~$9$,
\label{it:max9}
\item the vertex $u_1$ of $P_1\cap P_2$ has no red neighbour in $U'$, 
every vertex of $P_1-u_1$ has at most $5$ red neighbours in $U'$,
every vertex of $P_2-u_1$ has at most $3$ red neighbours in $U'$ except the neighbour of $u_1$ on $P_2$ which may have up to $4$ red neighbours in $U'$, and%
\label{it:P12red}
\item at the end of the partial sequence, that is for $U^0:=V(G^0)\setminus(V(C)\cup W)$, we have \mbox{$\big|U^0\cap L_{\ell+1}[G^0]\big|\leq2$}
and $\big|U^0\cap L_j[G^0]\big|\leq1$ for all $j>\ell+1$.
\label{it:contractfin}
\end{enumerate}
\end{lemma}

\begin{figure}[tb]
$$
\begin{tikzpicture}[scale=1]
\small
\tikzstyle{every node}=[draw, shape=circle, minimum size=3pt,inner sep=0pt, fill=black]
\node at (0,0) (x) {};
\node at (6,0) (y) {};
\node[label=above:$r$] at (3,5) (z) {};
\draw (x)--(y) to[bend right] (z) (z) to[bend right] (x);
\draw (x)-- +(0.4,0.3); \draw (x)-- +(0.2,0.5);
\draw (y)-- +(-0.4,0.3); \draw (y)-- +(-0.2,0.5);
\node[label=left:$v_1$] at (1.75,0.75) (v1) {};
\node[label=right:$v_2$] at (3.75,1) (v2) {};
\node[label=right:$~{u_1=u_2}$] at (2.75,4) (u1) {};
\draw[thick, fill=gray!20] (v2) to[bend right] (u1) to[bend right] (v1) -- (v2);
\draw (z)--(3.2,4.7) node{}--(2.5,4.3) node{}--(u1);
\draw (v1)-- +(-0.2,-0.3); \draw (v1)-- +(0.2,-0.3);
\draw (v2)-- +(-0.2,-0.3); \draw (v2)-- +(0.2,-0.3);
\draw (u1)-- +(0,-0.4);
\node at (2.35,3.5){}; \node at (3.23,3.5){};
\node at (2.0,3){}; \node at (3.57,3){};
\node[label=left:$P_1~$] at (1.8,2.5){}; \node[label=right:$~P_2$] at (3.77,2.5){};
\node at (1.69,2){}; \node at (3.85,2){};
\node[label=above:$C\quad$] at (v2) {};
\node[draw=none,fill=none] at (2.75,2) {$U$};
\node[draw=none,fill=none] at (2.75,0.7) {$f$};
\node[draw=none,fill=none] at (5,1.5) {$W$};
\end{tikzpicture}
\raise15ex\hbox{\huge\boldmath~$\leadsto$~}
\begin{tikzpicture}[scale=1]
\small
\tikzstyle{every node}=[draw, shape=circle, minimum size=3pt,inner sep=0pt, fill=black]
\node at (0,0) (x) {};
\node at (6,0) (y) {};
\node[label=above:$r$] at (3,5) (z) {};
\draw (x)--(y) to[bend right] (z) (z) to[bend right] (x);
\draw (x)-- +(0.4,0.3); \draw (x)-- +(0.2,0.5);
\draw (y)-- +(-0.4,0.3); \draw (y)-- +(-0.2,0.5);
\node[label=left:$v_1$] at (1.75,0.75) (v1) {};
\node[label=right:$v_2$] at (3.75,1) (v2) {};
\node[label=right:$~{u_1=u_2}$] at (2.75,4) (u1) {};
\draw[thick] (v2) to[bend right] (u1) to[bend right] (v1) -- (v2);
\draw (z)--(3.2,4.7) node{}--(2.5,4.3) node{}--(u1);
\draw (v1)-- +(-0.2,-0.3); \draw (v1)-- +(0.2,-0.3);
\draw (v2)-- +(-0.2,-0.3); \draw (v2)-- +(0.2,-0.3);
\node at (2.35,3.5) (pa){}; \node at (3.23,3.5) (qa){};
\node at (2.0,3) (pb){}; \node at (3.57,3) (qb){};
\node[label=left:$P_1~$] at (1.8,2.5) (pc){}; \node[label=right:$~P_2$] at (3.77,2.5) (qc){};
\node at (1.69,2) (pd){}; \node at (3.85,2) (qd){};
%
\node[draw=none,fill=none] at (2.7,1.5) {$U^0$};
\node[draw=none,fill=none] at (5,1.5) {$W$};
\node[draw=none,fill=none,red] at (3.3,1.75) {\scriptsize$(\leq\!9)$};
\node[draw=none,fill=none,red] at (1.5,3) {\scriptsize$(\leq\!5)$};
\node[draw=none,fill=none,red] at (3.66,3.5) {\scriptsize$(\leq\!4)$};
\node[draw=none,fill=none,red] at (4.05,3) {\scriptsize$(\leq\!3)$};
\node[red] at (2.7,3.5) (ra) {}; \node[red] at (2.93,3.5) (rb) {};
\node[red] at (2.8,3) (rc) {}; \node[red] at (2.6,2.5) (rd) {}; \node[red] at (2.83,2) (re) {};
\draw (u1)--(ra);
\draw[red, thick,densely dotted] (pa)--(qa) (pb)--(qb)  (pc)--(qc) (pd)--(qd) (re)--(rd)--(rc)--(rb) (rc)--(ra);
\draw[red, thick,densely dotted] (ra) to[bend left] (qa) (rb) to[bend left] (pa);
\draw[red, thick,densely dotted] (pc)--(rc)--(qa) (pd)--(rd)--(qb) (ra)--(pb)--(rb) (rc)--(pa) (rd)--(pb) (re)--(pc);
\end{tikzpicture}
$$
\caption{(left) The setup of \Cref{lem:core}, where $v_1$ is to the left of~$v_2$.
	(right) The outcome of the claimed partial contraction sequence which contracts only vertices of $U$ inside the shaded region
	from the left and maintains bounded red degrees in the region and on its boundary $C$ formed by $P_1$ and $P_2$ and $f$.
	Not all depicted red edges really exist, and some of them may actually be black.}
\label{fig:corelem}
\end{figure}
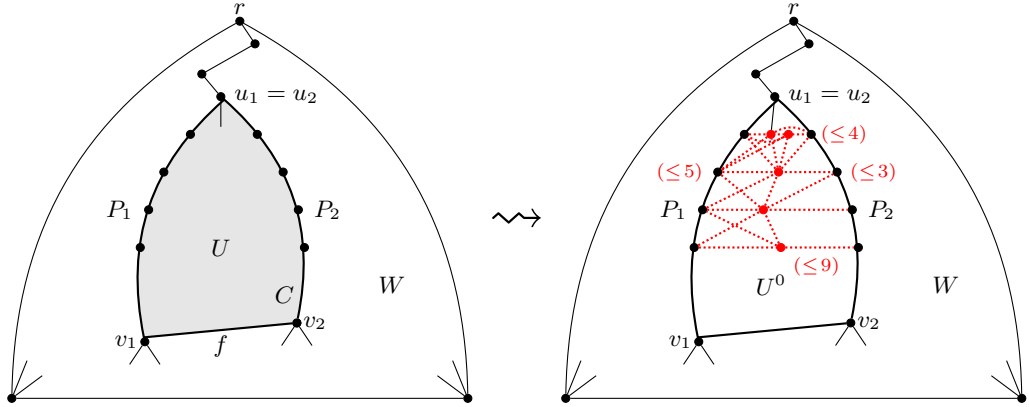

Before proceeding further, we comment on two important things.
First, the possible case of \mbox{$\big|U^0\cap L_{\ell+1}[G^0]\big|=2$} in \eqref{it:contractfin} of \Cref{lem:core} is unavoidable since
a contraction of this pair could result in a red edge from $u_1$ which must be forbidden in the induction.
Second, we remark that, in \Cref{lem:core}, all vertices of $U$ have distance from $r$ greater than~$\ell$, but
also the distance from $r$ of some vertices in $U$ may be greater than the distance from $r$ to $v_1$ or~$v_2$ (and our coming proof is aware of this possibility).
We illustrate \Cref{lem:core} in \Cref{fig:corelem}.

We also observe that the assumptions of \Cref{lem:core} directly imply some other properties useful for the proof of \Cref{lem:core}.

\begin{claim}\label{clm:redalign}\it
Within the notation and assumptions of \Cref{lem:core}, for $G'$ we also have that:
\begin{enumerate}[(I)]\setcounter{enumi}{3}
\item Every red edge in $G'$ has one end in $U'$ and the other end in~$U'\cup V(C)$,
\label{it:nooutside}
\item if $v\in V(P_2)\cap L_j[G']$, then there is no edge in $G'$ (red or black) from $v$ to a vertex of $(U'\cup V(P_1))\cap L_{j-1}[G']$.
\label{it:redalign2}
\end{enumerate}
\end{claim}
\begin{claimproof}
Regarding \eqref{it:nooutside}, observe that since only vertices that stem from $U$ participate in contractions, every red edge of $G'$ must have an end in $U'$.
Furthermore, due to the cycle $C$ in the plane graph $G$, no vertex of $W$ is adjacent to a vertex of $U$ in $G$, and so the same holds in $G$.
Then no vertex of $W$ is ever adjacent to a vertex being contracted in our sequence from $G$ to $G^0$ (which includes~$G'$).

Concerning \eqref{it:redalign2}, if $v$ was adjacent to $x'\in U'\cap L_{j-1}[G']$, then, among the vertices of $U\cap L_{j-1}$ contracted into~$x'$,
there had to be a neighbour $x\in U\cap L_{j-1}$ of $v$ in $G$, and so also $xv\in E(G)$.
Then $x$ was to the left of $v$ in $G$ by the definition of $U$, which contradicts the assumption that $T$ is left-aligned.
\end{claimproof}

\begin{claim}\label{clm:redbfin}\it
With the notation and assumptions of \Cref{lem:core}, we have that if \eqref{it:contractfin} is true, 
then \eqref{it:P12red} and \eqref{it:max9} except the condition on $u_1$ are also true for~$G'=G^0$.
\end{claim}
\begin{claimproof}
Every vertex of $P_1-u_1$ may have red neighbours only in $U^0$ by \eqref{it:nooutside}, and by \eqref{it:contractfin} and \Cref{cl:layeri3} there are at most $2+1+1=4$ 
such potential neighbours (layer by layer, and this can be improved further).
Likewise, every vertex of $P_2-u_1$ may have red neighbours only in $U^0$, and using additionally \eqref{it:redalign2}, there are at most $2+1=3$ such potential neighbours.
For every vertex of~$U^0$, we have at most $4+2+3=9$ such potential neighbours simply by \eqref{it:contractfin} and \Cref{cl:layeri3}
(again layer by layer, and the bound can be improved further).
\end{claimproof}

We also show how \Cref{lem:core} implies the first part of our main result:

\begin{proof}[Proof of \Cref{thm:twwplanar} (the upper bound)]
We start with a given simple planar graph $H$, and extend any plane embedding of $H$ into a simple plane triangulation $G$
such that $H$ is an induced subgraph of~$G$.
Then we choose a root $r$ on the outer face of $G$ and, for some left-aligned BFS tree of $G$ rooted at $r$ which exists by \Cref{clm:existslal},
and the paths $P_1$ and $P_2$ formed by each of the two edges of the outer face incident to $r$, apply \Cref{lem:core}.
This way we get a partial contraction sequence from $G$ to a trigraph $G^0$ of maximum red degree $9$ (along the sequence).
Observe that $L_0[G^0]=\{r\}$ and $L_1[G^0]$ contains the two vertices of the outer face other than $r$ and at most 
two vertices of $U^0$ by \eqref{it:contractfin}.
For all $j\geq2$, we have $\big|L_j[G^0]\big|\leq1$ again by \eqref{it:contractfin}.

In the final phase, we pairwise contract the at most $4$ vertices of $L_1[G^0]$ into one vertex, making $G^0$ into a path by \Cref{cl:layeri3},
and then naturally contract this path down to a vertex, never exceeding red degree $5<9$ along this final phase.
The restriction of the whole contraction sequence of $G$ to $V(H)$ then certifies that the twin-width of $H$ is at most~$9$, too.
\end{proof}

Now we get to the proof of core \Cref{lem:core} which will conclude our main result.

\begin{proof}[Proof of \Cref{lem:core}]
We first resolve several possible degenerate cases.
If $U=\emptyset$, we are immediately done with the empty partial contraction sequence. So, assume $U\not=\emptyset$.

Recall the edge $f=\{v_1,v_2\}\in E(G)$ connecting the other ends of $P_1$ and~$P_2$.
If $v_1$ has no neighbour in $U$, then $\{v_2,v_3\}\in E(G)$ where $v_3$ is the neighbour of $v_1$ on $P_1$.
In such case, we simply apply \Cref{lem:core} inductively to $P_1-v_1$ and~$P_2$, while the rest of the assumptions remain the same.
The symmetric argument is applied when $v_2$ has no neighbour in $U$.

Otherwise, let $v_3\in U$ be the vertex (unique in $U$) such that $(v_1,v_2,v_3)$ bound a triangular face of $G$.
Let $P\subseteq T$ be the vertical path connecting $v_3$ to the root $r$ where, by $r\not\in U$ and planarity of $G$, we have $u_1\in V(P)$.
If $v_1\in V(P)$, then we (similarly as above) apply \Cref{lem:core} inductively to $P$ and $P_2$, while the rest of the assumptions again remain the same.
In the resulting trigraph $G^1$, assuming $v_3\in L_j$ (note that $j\geq\ell+2$) and $L_j\cap U\not=\{v_3\}$ at the beginning,
we then contract $v_3$ with the (unique by \eqref{it:contractfin}) vertex of $L_j[G^1]\cap U^1$ where $U^1:=V(G^1)\setminus(V(P\cup P_2)\cup W)$.
This is correct by \Cref{clm:redbfin}.

In all other cases we continue as follows. See in \Cref{fig:corerec}.
\medskip

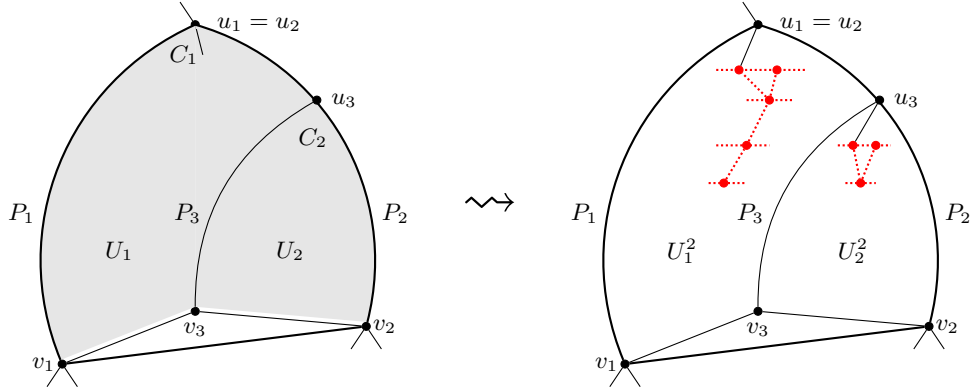
\begin{figure}[th]
$$
\begin{tikzpicture}[scale=1]
\small
\tikzstyle{every node}=[draw, shape=circle, minimum size=3pt,inner sep=0pt, fill=black]
\node[label=left:$v_1$] at (1,0.5) (v1) {};
\node[label=right:$v_2$] at (5,1) (v2) {};
\node[label=right:$~~{u_1=u_2}$] at (2.75,5) (u1) {};
\node[label=below:$v_3$] at (2.75,1.2) (v3) {};
\draw[draw=none, fill=gray!20] (v1)--(v3)--(2.75,5) to[bend right=44] (v1);
\draw[draw=none, fill=gray!20] (v3)--(v2) to[bend right=44] (2.75,5)--(v3);
\draw[thick] (v2) to[bend right=44] (2.75,5) to[bend right=44] (v1) (v1) -- (v2);
\draw (v2)--(v3)--(v1);
\node[label=right:$~u_3$] at (4.35,4) (u3) {};
\draw (v3) to[bend left] (u3);
\draw (u1)-- +(-0.2,0.3) (u1)-- +(0.1,-0.4);
\draw (v1)-- +(-0.2,-0.3); \draw (v1)-- +(0.2,-0.3);
\draw (v2)-- +(-0.2,-0.3); \draw (v2)-- +(0.2,-0.3);
\node[draw=none,fill=none, label=left:$P_1~$] at (0.85,2.5){};
\node[draw=none,fill=none, label=right:$~P_2$] at (5,2.5){};
\node[draw=none,fill=none] at (2.65,2.5) {$P_3$};
\node[draw=none,fill=none] at (1.75,2) {$U_1$};
\node[draw=none,fill=none] at (4,2) {$U_2$};
\node[draw=none,fill=none] at (4.3,3.5) {$C_2$};
\node[draw=none,fill=none] at (2.6,4.6) {$C_1$};
\end{tikzpicture}
\raise15ex\hbox{\huge\boldmath\quad$\leadsto$\quad}
\begin{tikzpicture}[scale=1]
\small
\tikzstyle{every node}=[draw, shape=circle, minimum size=3pt,inner sep=0pt, fill=black]
\node[label=left:$v_1$] at (1,0.5) (v1) {};
\node[label=right:$v_2$] at (5,1) (v2) {};
\node[label=right:$~~{u_1=u_2}$] at (2.75,5) (u1) {};
\node[label=below:$v_3$] at (2.75,1.2) (v3) {};
\draw[thick] (v2) to[bend right=44] (2.75,5) to[bend right=44] (v1) (v1) -- (v2);
\draw (v2)--(v3)--(v1);
\node[label=right:$~u_3$] at (4.35,4) (u3) {};
\draw (v3) to[bend left] (u3);
\draw (u1)-- +(-0.2,0.3);
\draw (v1)-- +(-0.2,-0.3); \draw (v1)-- +(0.2,-0.3);
\draw (v2)-- +(-0.2,-0.3); \draw (v2)-- +(0.2,-0.3);
\node[draw=none,fill=none, label=left:$P_1~$] at (0.85,2.5){};
\node[draw=none,fill=none, label=right:$~P_2$] at (5,2.5){};
\node[draw=none,fill=none] at (2.65,2.5) {$P_3$};
\node[draw=none,fill=none] at (1.75,2) {$U^2_1$};
\node[draw=none,fill=none] at (4,2) {$U^2_2$};
\node[red] at (2.5,4.4) (ra) {}; \node[red] at (3,4.4) (rb) {};
\node[red] at (2.9,4) (rc) {}; \node[red] at (2.6,3.4) (rd) {}; \node[red] at (2.3,2.9) (re) {};
\draw (u1)--(ra);
\node[red] at (4,3.4) (sa) {}; \node[red] at (4.3,3.4) (sb) {};
\node[red] at (4.1,2.9) (sc) {};
\draw (u3)--(sa);
\draw[red, thick,densely dotted] (re)--(rd)--(rc)--(rb) (rc)--(ra) (sa)--(sc)--(sb);
\draw[red, thick,densely dotted] (2.2,4.4)--(3.4,4.4) (2.6,4)--(3.2,4) (2.2,3.4)--(2.9,3.4) (2.1,2.9)--(2.6,2.9);
\draw[red, thick,densely dotted] (3.8,3.4)--(4.5,3.4) (3.9,2.9)--(4.3,2.9);
\end{tikzpicture}
$$
\caption{(left) The setup of the main induction step in the proof of \Cref{lem:core} (another possibilities are that $u_3$ lies on $P_1$, or is equal to~$u_1$).
	(right) A schematic outcome of the recursive application of \Cref{lem:core} to each of the cycles $C_1$ and $C_2$.
	Note that the only vertices that have neighbours in both contracted subsets $U_1^2$ and $U_2^2$ are those of $P_3$,
	and among them only $u_3$ has possibly neighbours not explicitly mentioned in this picture, that is, neighbours not in
	$V(P_1\cup P_2\cup P_3)$ and not among the vertices that stem from $U_1\cup U_2$.}
\label{fig:corerec}
\end{figure}

Observe that $v_3$ is to the left of $v_2$, and $v_1$ is to the left of $v_3$.
Let $P_3\subseteq P$ be the subpath of $P$ from $v_3$ to the first vertex $u_3\in V(P)\cap V(C)$ shared with the cycle~$C$.
Assuming $u_3\in V(P_k)$ for $k\in\{1,2\}$ (and the following is sound even if $u_3=u_1\in V(P_1)\cap V(P_2)$\,),
let $P_k^-$ be the subpath of $P_k$ from $v_k$ to $u_3$,
and $P_3^+\subseteq P$ be the subpath of $P$ from $v_3$ to $u_1$ (which is the union of $P_3$ with the subpath of $P_k$ from $u_3$ to~$u_1$).

We have got two cycles $C_1,C_2\subseteq G$, namely $C_k:=(P_k^-\cup P_3)+\{v_k,v_3\}$ and $C_{3-k}:=(P_{3-k}\cup P_3^+)+\{v_{3-k},v_3\}$.
For $i=1,2$, we denote by $U_i\subseteq U$ the set of vertices drawn in the interior of $C_i$.
Note that $U$ is tri-partitioned into $U_1$, $U_2$ and $V(P_3)\setminus\{u_3\}$.
We apply \Cref{lem:core} inductively to each of the cycles $C_1$ and $C_2$. This is sound by the following:
\begin{itemize}
\item Each of the paths $P_3$, $P_3^+$ and $P_k^-$ is vertical wrt.~$T$ (as a subpath of other vertical path).
Furthermore, $V(P_k^-)\cap V(P_3)=\{u_3\}$ and $V(P_{3-k})\cap V(P_3^+)=\{u_1\}$.
\item We get that each of the vertical paths $P_k^-$, $P_3$ and $P_3^+$ is of length at least~$1$,
and it is possible to apply \Cref{lem:core} to both $C_1$ and $C_2$ independently and concurrently.
\item Since $v_3\in U$ belongs to neither of $U_1,U_2$, we always have $|U_1|+|U_2|<|U|$.
\end{itemize}

This way we get $\ca L$-respecting partial contraction sequences $\tau_1$ (for $C_1$ contracting in $U_1$) 
and $\tau_2$ (for $C_2$ contracting in $U_2$) of $G$.
These sequences are independent in the sense that no contraction pair in $U_1$ creates a red edge
incident to a vertex of $U_2$ by inductive invocation of condition \eqref{it:nooutside}, and vice versa.
We may thus append $\tau_2$ after~$\tau_1$, as the concatenation $\tau:=\tau_1.\tau_2$, such that
$\tau$ is a valid $\ca L$-respecting partial contraction sequence of $G$ resulting in a trigraph $G^2$.
At every step of $\tau$, each vertex that stems from $U_1$ or $U_2$ is of red degree at most $9$ by inductive invocation of condition \eqref{it:max9}.
Furthermore, by \eqref{it:P12red}, each vertex of $V(P_3)\setminus\{u_3\}$ is of red degree at most $5+4=9$ at every step of $\tau$.
This fulfills \eqref{it:max9} for~$\tau$.

Concerning \eqref{it:P12red}, the property is inherited for $\tau$ from its recursive application to one of $\tau_1$ or $\tau_2$, using also \eqref{it:nooutside},
for every vertex of $P_1\cup P_2$ (including for $u_1$), except possibly for $u_3$ which belongs to both bounding cycles $C_1$ and $C_2$.
However, condition \eqref{it:P12red} of $u_3$ from the recursive application to $\tau_k$ claims that there are no red edges created from $u_3$ to vertices that stem from $U_k$,
and so \eqref{it:P12red} is inherited also for the vertex $u_3$ from its application to~$\tau_{3-k}$.
Condition \eqref{it:contractfin} is void for this intermediate sequence~$\tau$.

\medskip
We are going to continue with $\ca L$-respecting partial contraction sequences $(\tau_0.\tau_3)$ from $G^2$ to
$G^3$, then $\tau_4$ from $G^3$ to $G^4$ and (possibly empty) $\tau_5$ from $G^4$ to $G^5$, such that 
the concatenation $(\tau.\tau_0.\tau_3.\tau_4.\tau_5)$ fulfills all conditions of \Cref{lem:core} for the cycle~$C$.
For $i=1,2$, denote by $U_i^2\subseteq V(G^2)$ the set of vertices which stem from $U_i$ by the sequence~$\tau$.
Recall that $u_3\in V(P_k)$ where $k\in\{1,2\}$, and let $a>\ell$ be such that $u_3\in L_{a-1}$.
If $u_3=u_1$ or $|L_{a}[G^2]\cap U_k^2|\leq1$, then let $\tau_0=\emptyset$.
Otherwise, $a\geq\ell+2$ and the partial contraction sequence~$\tau_0$ has one step contracting the vertex pair $L_{a}[G^2]\cap U_k^2$.

The contraction in $\tau_0$ is safe under the stated conditions, as we now show.
First, if $k=1$, then $u_3$ actually has red degree at most $2+1+1=4$ in $G^2$ 
(counting BFS layers $L_{a-2}[G^2]$, $L_{a-1}[G^2]$ and $L_a[G^2]$ among the vertices of~$U_2^2$),
and so $u_3$ may accept a new red neighbour from this contraction while keeping \eqref{it:P12red} true.
The red degrees of all other affected vertices, including the contracted one, remain good in $\tau_0$ by \Cref{clm:redbfin} (as applied to~$C_1$).
If $k=2$ and $a\geq\ell+3$, then the situation is again good since, in $G^2$ before the $\tau_0$ contraction, the vertex $u_3\in L_{a-1}$ has at most $1+1=2$
red neighbours in $\big(L_{a}[G^2]\cup L_{a-1}[G^2]\big)\cap U_1^2$ and no one in $L_{a-2}[G^2]$ by \eqref{it:redalign2} (as applied to the subcase of $C_1$), 
and so $u_3$ may accept a new red neighbour from the $\tau_0$ contraction while keeping \eqref{it:P12red} true.
If $a=\ell+2$, then we actually have $u_3\in L_{\ell+1}$, i.e., $u_3$ is the neighbour of $u_1$ on~$P_2$.
Then $u_3$ has at most $3$ red neighbours which all belong to $\big(L_{\ell+2}[G^2]\cup L_{\ell+1}[G^2]\big)\cap U_1^2$, 
and by the special provision in \eqref{it:P12red} $u_3$ may accept a fourth red neighbour from the $\tau_0$ contraction, too.

Subsequently, the partial contraction sequence~$\tau_3$ contracts vertex pairs from $U_2^2\cup V(P_3)\setminus\{u_3\}$ (but it possibly skips some vertices which are close to $u_1$),
the partial contraction sequence $\tau_4$ contracts vertex pairs from $U_1^2\cup V(P_3)\setminus\{u_3\}$ (again possibly skipping some vertices which are close to $u_1$),
and the sequence $\tau_5$ then contracts pairs among the possibly skipped vertices from $\tau_3$ and~$\tau_4$.

We continue with the definition of $\tau_3$.
Let $d:=\max\{i: (U_2^2\cup\{v_3\})\cap L_{i}[G^2]\not=\emptyset\}$.
Let $G^2_0$ be the trigraph resulting from $\tau_0$ applied to~$G^2$.
Note that, for $a\leq j\leq d$ we have $|L_{j}[G^2_0]\cap V(P_3)|\leq1$, and for $\ell+2\leq j\leq d$ we have $|L_{j}[G^2_0]\cap U_1^2|\leq1$ and $|L_{j}[G^2_0]\cap U_2^2|\leq1$.
%
%
For $j=d,d-1,\ldots,b=\max(\ell+3,a)$ in this order, $\tau_3$ contracts the pair $L_{j}[G^2_0]\cap(U_2^2\cup V(P_3))$,
or $\tau_3$ skips the step if this set has at most~$1$ vertex. That is all for $\tau_3$.

Let $G''$ be an (intermediate) trigraph along the sequence $\tau_3$, and $U'':=V(G'')\setminus(V(C)\cup W)$.
We now verify validity of conditions \eqref{it:P12red} and \eqref{it:max9} for~$G''$.
This is trivial for a vertex $x\in V(P_2)$ since no new potential red neighbour has been created along $\tau_3$.
On the other hand, the case is more complex for a vertex $x\in V(P_1)$ since such $x$ may gain new red neighbours which stem from $V(P_3)\cup U_2^2$.
However, assuming $x\in L_c$ and $x\not=u_3$, there is no edge in $G$ from $L_{c+1}\cap V(P_3)$ to $x$ by~\eqref{it:redalign2},
and also no such edge in $G^2_0$ from $L_{c+1}[G^2_0]\cap U_2^2$ by planarity, and consequently, 
no new red edge to $x$ is created by the $\tau_3$-contraction in layer $c+1$.
If $x=u_3$, then the edge of $P_3$ from $x=u_3$ turns red by the $\tau_3$-contraction in layer $c+1$, but still $x$ has at most
$3$ red neighbours in $(L_{c-1}[G^2_0]\cup L_c[G^2_0])\cap U_2^2$ and at most $2$ in layer $c+1$ after this contraction.
Hence the situation is good if $\tau_3$ has not yet contracted in layer $c$ in~$G''$, and if a contraction in layer $c$ has happened, then we know that~$c\geq\ell+3$.
In the latter case, in $G''$ there are up to $3$ red neighbours of $x$ in $U_1^2$, up to $2$ among the vertices that stem from $V(P_3)\cup U_2^2$
and no more, as needed in \eqref{it:P12red}.

Among vertices $x\in U''$, the case is again trivial for $x\in U_1^2$ since no new potential red neighbour of such $x$ has been created along $\tau_3$.
So, consider $x$ which is in or stems from $V(P_3)\cup U_2^2$, and assume that~$x\in L_c[G'']$.
Then, no matter of where $G''$ is along the sequence, $x$ has at most $3$ red neighbours in layer $c-1$ since one can be in $V(P_2)\cap L_{c-1}[G'']$
and at most two in $(U''\setminus U_1^2)\cap L_{c-1}[G'']$. 
There cannot be any neighbour of $x$ in $(V(P_1)\cup U_1^2)\cap L_{c-1}[G'']$ because of planarity and \eqref{it:redalign2}.
Likewise, $x$ has at most $3$ red neighbours in layer $c$ simply since there are no more candidates in the layer.
Regarding layer $c+1$, again by \eqref{it:redalign2}, there cannot be a red edge from $x$ to $V(P_2)\cap L_{c+1}[G'']$.
So, before $\tau_3$ contracts in layer~$c+1$, $x$~has at most two red neighbours there, and after that at most three such
(considering also a possible $\tau_3$-contraction in layer $c$).
This sums to at most $9$, as needed in \eqref{it:max9}.

We are done with the partial sequence $\tau_3$.

\medskip
We continue with the definition of the partial contraction sequence~$\tau_4$ starting from the trigraph $G^3$ which is the result of~$\tau_3$.
Denote by $U_2^3$ the vertices of $G^3$ that stem from $V(P_3)\cup U_2^2$ along $\tau_3$, and note that for $b\leq j\leq d$ we have $|L_{j}[G^3]\cap U_2^3|\leq1$.
We intend to proceed similarly to the sequence $\tau_3$ (but with some changes caused by the fact that the paths $P_1$ and $P_2$ are not symmetric in our assumptions).
Let analogously $d':=\max\{i: (U_1^2\cup\{v_3\})\cap L_{i}[G^3]\not=\emptyset\}$.
%
%
For $j=d',d'-1,\ldots,b=\max(\ell+3,a)$ in this order, $\tau_4$ contracts the pair \mbox{$L_{j}[G^3]\cap(U_1^2\cup U_2^3)$},
or skips the step if this set has at most $1$ vertex. This finishes $\tau_4$.

Let $G''$ now be an (intermediate) trigraph along the sequence $\tau_4$, and $U'':=V(G'')\setminus(V(C)\cup W)$.
We again verify validity of conditions \eqref{it:P12red} and \eqref{it:max9} for~$G''$.
This time the check is trivial for vertices $x\in V(P_1)\cup V(P_2)$ since no new potential red neighbour
(i.e., not counted in our previous arguments along $\tau_3$) for them has been created along $\tau_4$.
If $x\in U''$ such that $x\in L_{c}[G'']$ and $c\leq b-1$, then the argument of no new potential red neighbour also applies along~$\tau_4$.

Assume $x\in U''$ such that $x\in L_{c}[G'']$ and $c>b$.
If the partial contraction sequence $\tau_4$ has not yet reached layer $c+1$, then the red neighbours of $x$ in $G''$ are as in $G^3$,
and so we assume that the $\tau_4$-contractions have already reached layer $c+1$.
Hence, $x$ has at most $2$ red neighbours in layer $c+1$ by \eqref{it:redalign2}, 
and since $\big|(V(C)\cup U'')\cap L_{i}[G'']\big|\leq4$ for $i\geq c-1\geq b$ (even if $\tau_4$-contractions have not yet reached layer~$c$), 
$x$ has at most $3+4=7$ red neighbours in layers $c$ and $c-1$, together at most~$9$ as desired.
Special care has to be taken of the case $c=b$, in which $x$ may potentially have up to $5$ red neighbours in layer $b-1\geq\ell+2$.
If $\tau_4$ is not at the end, i.e.~the layer $c=b$ has not been contracted yet,
then $x\in U^2_1\cap L_{c}[G'']$ has no neighbour on $P_2$ in layer $c-1$ by planarity,
and $x\in U^3_2\cap L_{c}[G'']$ has no neighbour on $P_1$ in layer $c-1$ by planarity and \eqref{it:redalign2} as applied during the sequence~$\tau_3$.
In both possibilities of $x$ we thus get an upper bound of $2+3+4=9$.
At the end of $\tau_4$, $x$ has no other potential red neighbour in layer $c$ than the two on $P_1$ and $P_2$, and we hence get~$2+2+5=9$.

We are finished with the partial contraction sequence~$\tau_4$.

\medskip
To recapitulate, we have verified conditions \eqref{it:P12red} and \eqref{it:max9} for the partial contraction sequence $(\tau.\tau_0.\tau_3.\tau_4)$,
and we have also got that condition \eqref{it:contractfin} is true for all~$j\geq b=\max(\ell+3,a)$.
Furthermore, \eqref{it:contractfin} is true for all~$j<a$ by the initial inductive application of \Cref{lem:core} to one of $C_1$ or $C_2$ (precisely, to~$C_{3-k}$).
So, we are completely done (with empty~$\tau_5$) unless $a<b$, where the latter can happen (recall our choice of $b$) only if $b\leq\ell+3$ and~$a\leq\ell+2$.
In other words, either when $u_3=u_1$ or when $u_3$ is the neighbour of $u_1$ on $P_1$ or on~$P_2$.
A straightforward case check reveals that there are only three such cases depicted in \Cref{fig:endcases}.
If none of these happens, we simply set $\tau_5=\emptyset$.

In the case (a) of \Cref{fig:endcases}, $\tau_5$ contracts $x_3$ with $x_4$, and then with~$z$, which does not exceed red degrees $9$,
specifically for the vertex which stems from $x_3$ thanks to \eqref{it:redalign2}.
In the case (b), $\tau_5$ again contracts $x_3$ with $x_4$, then with~$z$, and the argument is the same.

Situation is different in the case (c) of $u_1=u_3$, in which we may up to symmetry assume that $x_1u_1,x_4u_1\in E(G)$ and $x_2u_1,x_3u_1\not\in E(G)$.
Then $\tau_5$ in the first two steps contracts the pair $x_1,x_4$ and the pair $x_2x_3$ (if these pairs exist).
This does not create red edges from $u_1$, and in particular, the red degree of the vertex obtained by the contraction of $x_1,x_4$
stays at most $9$ ($5$ in layer $\ell+1$ and $4$ in layer $\ell+2$) thanks to $x_4$ not having a red neighbour on $P_2$ in layer $\ell+2$ by \eqref{it:redalign2}.
The red degrees of other concerned vertices, including those on $P_1$ and $P_2$, are easily checked from the picture.
Next to that, $\tau_5$ contracts $t$ with the vertex that stems above from $x_1,x_4$ if $\{t,u_1\}\in E(G)$, 
or $\tau_5$ contracts $t$ with the vertex that stems from $x_2,x_3$ if $\{t,u_1\}\not\in E(G)$, and this again does not create a red edge towards $u_1$
while keeping other concerned red degrees within the required bounds.
We finish $\tau_5$ by contracting $x_5$ with $x_6$, and then with $z$ as in the cases (a) and (b) and using the same argument.

\begin{figure}[th]
$$\mbox{a)\hspace*{-2ex}}
\begin{tikzpicture}[scale=1]
\small
\tikzstyle{every node}=[draw, shape=circle, minimum size=3pt,inner sep=0pt, fill=black]
\node[draw=none,fill=none, label=below:$P_1$] at (1,0.9) (v1) {};
\node[draw=none,fill=none, label=below:$P_2$] at (5,1) (v2) {};
\node[label=right:$~~{u_1}$] at (2.75,5) (u1) {};
\node[label=left:$u_3~$] at (1.45,4) (u3) {};
\node[draw=none,fill=none] at (3,3) (v3) {};
\draw[thick] (v2) to[bend right=44] (2.75,5) to[bend right=44] (v1);
\node at (4.35,4) (r1) {}; \node at (0.92,3) (p2) {}; \node[fill=none] at (3,3) (q2) {\,$z$\,}; \node at (4.9,3) (r2) {};
\draw[thick] (q2) to[bend right=10] (u3);
\node[draw=none,fill=none] at (0,4) {($\ell+1$)};
\node[draw=none,fill=none] at (0,3) {($\ell+2$)};
\node[draw=none,fill=none] at (0,2) {($\ell+3$)};
\node[draw=none,fill=none] at (1.75,1.5) {($U^2_1$)};
\node[draw=none,fill=none] at (4,1.5) {($U^3_2$)};
\tikzstyle{every node}=[draw, color=red, shape=circle, minimum size=12pt,inner sep=0pt, fill=none]
\node at (2.5,4) (ra) {$x_1$};  \node at (3.5,4) (rb) {$x_2$};
\node at (2.1,3) (rd) {$x_3$};  \node at (4.1,3) (rc) {$x_4$};
\node at (3.3,2) (re) {$x_5$};
\draw (u1)--(ra);
\tikzstyle{every path}=[draw, color=red, thick,densely dotted]
\draw (rb) to[bend right] (u3)--(ra)--(rb)--(r1) to[bend right] (ra);
\draw (u3)--(rc)--(ra)--(q2)--(rd)--(0.8,2) (q2)--(rb)--(rc)--(q2)--(rd)--(p2) (r1)--(rc);
\draw (p2)--(re)--(q2) (rc)--(re)--(rd)--(u3) (re)--(r2)--(rc) (0.8,2)--(re)--(5.1,2);
\draw (re)-- +(-0.5,-0.25) (re)-- +(0,-0.5);
\end{tikzpicture}
\mbox{\quad b)}
\begin{tikzpicture}[scale=1]
\small
\tikzstyle{every node}=[draw, shape=circle, minimum size=3pt,inner sep=0pt, fill=black]
\node[draw=none,fill=none, label=below:$P_1$] at (1,0.9) (v1) {};
\node[draw=none,fill=none, label=below:$P_2$] at (5,1) (v2) {};
\node[label=right:$~~{u_1}$] at (2.75,5) (u1) {};
\node[label=right:$~u_3$] at (4.35,4) (u3) {};
\node[draw=none,fill=none] at (3,3) (v3) {};
\draw[thick] (v2) to[bend right=44] (2.75,5) to[bend right=44] (v1);
\node at (1.45,4) (p1) {}; \node at (0.92,3) (p2) {}; \node[fill=none] at (3,3) (q2) {\,$z$\,}; \node at (4.9,3) (r2) {};
\draw[thick] (q2) to[bend left=10] (u3);
\node[draw=none,fill=none] at (1.75,1.5) {($U^2_1$)};
\node[draw=none,fill=none] at (4,1.5) {($U^3_2$)};
\tikzstyle{every node}=[draw, color=red, shape=circle, minimum size=12pt,inner sep=0pt, fill=none]
\node at (2.3,4) (ra) {$x_1$};  \node at (3.3,4) (rb) {$x_2$};
\node at (2.1,3) (rc) {$x_3$};  \node at (4.1,3) (rd) {$x_4$};
\node at (2.7,2) (re) {$x_5$};
\draw (u1)--(ra);
\tikzstyle{every path}=[draw, color=red, thick,densely dotted]
\draw (p2)--(rb) to[bend right] (p1)--(ra)--(rb)--(u3) to[bend right] (ra) (u3)--(rd);
\draw (u3)--(rc)--(ra)--(p2)--(rc)--(p1) (rb)--(rc)--(q2)--(rd)--(r2) (rc)--(0.8,2);
\draw (p2)--(re)--(q2) (rc)--(re)--(rd) (re)--(r2) (0.8,2)--(re)--(5.1,2);
\draw (re)-- +(-0.5,-0.25) (re)-- +(0,-0.5);
\end{tikzpicture}
$$ $$
\mbox{c)\hspace*{-2ex}}
\begin{tikzpicture}[xscale=1.5]
\small
\tikzstyle{every node}=[draw, shape=circle, minimum size=3pt,inner sep=0pt, fill=black]
\node[draw=none,fill=none, label=below:$P_1$] at (1,0.9) (v1) {};
\node[draw=none,fill=none, label=below:$P_2$] at (5,1) (v2) {};
\node[label=right:$~~{u_1=u_3}$] at (2.75,5) (u1) {};
\node at (2.75,5) (u3) {};
\node[draw=none,fill=none] at (3,3) (v3) {};
\draw[thick] (v2) to[bend right=44] (2.75,5) to[bend right=44] (v1);
\node at (1.44,4) (p1) {}; \node at (4.35,4) (r1) {}; \node at (0.92,3) (p2) {}; 
\node[fill=none] at (2.875,4) (q1) {\,$t$\,}; \node[fill=none] at (3,3) (q2) {\,$z$\,}; \node at (4.9,3) (r2) {};
\draw[thick] (q2)--(q1)--(u3);
\node[draw=none,fill=none] at (0.2,4) {($\ell+1$)};
\node[draw=none,fill=none] at (0.2,3) {($\ell+2$)};
\node[draw=none,fill=none] at (0.2,2) {($\ell+3$)};
\node[draw=none,fill=none] at (1.75,1.5) {($U^2_1$)};
\node[draw=none,fill=none] at (4,1.5) {($U^3_2$)};
\tikzstyle{every node}=[draw, color=red, shape=circle, minimum size=11pt,inner sep=0pt, fill=none]
\node at (1.85,4) (ra) {$x_1$};  \node at (2.4,4) (rb) {$x_2$};
\node at (3.35,4) (rc) {$x_3$};  \node at (3.9,4) (rd) {$x_4$};
\tikzstyle{every node}=[draw, color=red, shape=circle, minimum size=12pt,inner sep=0pt, fill=none]
\node at (2.1,3) (re) {$x_5$};  \node at (3.9,3) (rf) {$x_6$};
\node at (3.3,2) (rg) {$x_7$};
\draw (u1)--(ra) (u3)--(rd);
\tikzstyle{every path}=[draw, color=red, thick,densely dotted]
\draw (rb)--(ra)--(p1) to[bend left=45] (rb) (rb)--(q1) to[bend right=45] (ra);
\draw (rd)--(rc)--(q1) to[bend left=45] (rd) (rd)--(r1) to[bend right=45] (rc);
\draw (p2)--(ra)--(re)--(rb)--(p2)--(re)--(q2) (p1)--(re)--(q1);
\draw (q2)--(rc)--(rf)--(rd)--(q2)--(rf)--(r2) (q1)--(rf)--(r1);
\draw (p2)--(rg)--(q2) (re)--(rg)--(rf) (rg)--(r2) (0.8,2)--(rg)--(5.1,2);
\draw (re)--(0.8,2) (rg)-- +(-0.5,-0.25) (rg)-- +(0,-0.5);
\end{tikzpicture}
$$
\caption{The three possible cases at the end of the sequence $\tau_4$ in the proof of \Cref{lem:core}.
	In each case, all possible vertices resulting from the recursive contractions in layers $\ell+1$, $\ell+2$ and~$\ell+3$
	and all their potential red edges (dotted) are shown, but the three possibilities cover also (sub)cases in which some of 
	the red vertices or edges are missing, or edges are not red.
	For red edges not present in the pictures, consult \Cref{cl:layeri3} and \Cref{clm:redalign}.}
\label{fig:endcases}
\end{figure}
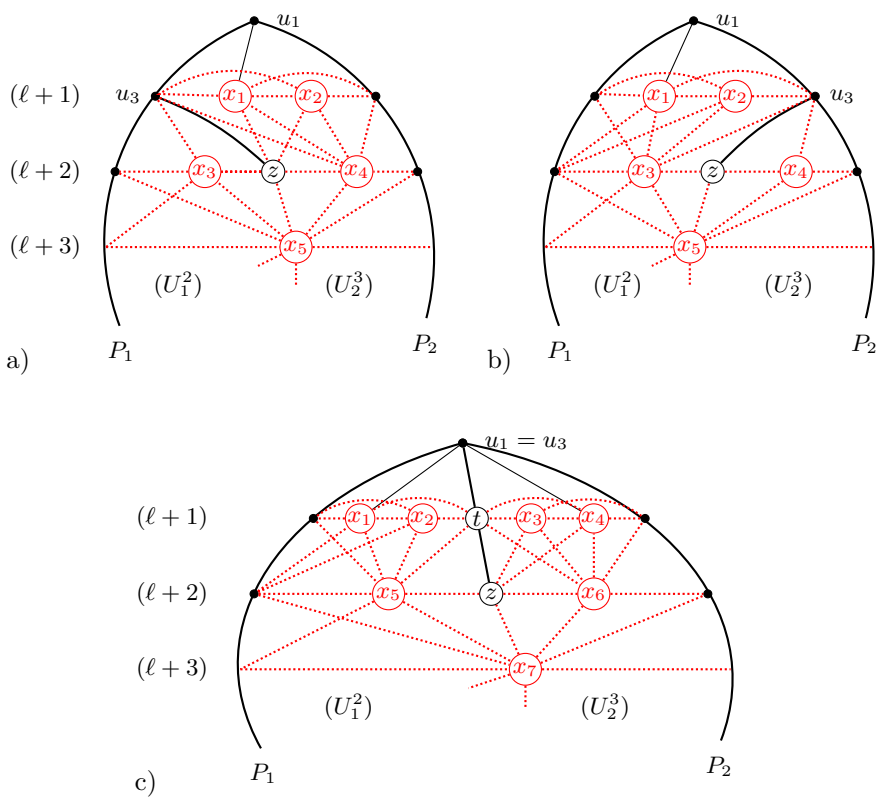

We have finally verified all conditions of \Cref{lem:core} for the concatenated sequence $(\tau.\tau_0.\tau_3.\tau_4.\tau_5)$, and so the lemma is proved.
\end{proof}

\begin{proof}[Proof of \Cref{thm:twwplanar} (the algorithmic part)]
We can construct a simple plane triangulation $G\supseteq H$ in linear time using standard planarity algorithms,
and then construct a left-aligned BFS tree $T\subseteq G$ again in linear time by \Cref{clm:existslal}.
In the rest, we straightforwardly implement the recursive division of $G$ as used in the proof of \Cref{lem:core},
and construct the contraction sequence of $H$ on return from the recursive calls as defined in the proof.
Note that we do not need at all to construct the intermediate trigraphs along the constructed contraction sequence,
and so the construction of the sequence is very easy---each recursive call returns just a simple list of the vertices
which stem from the recursive contractions, indexed by the BFS layers.
Then these (up to) two lists are easily in linear time ``merged'' together with the dividing path $P_3$, 
as specified by the proof of \Cref{lem:core}, into the resulting list of this call.

We may account total runtime in the ``division part'' of the algorithm to the edge(s) of $v_3$ into $v_1$ or $v_2$ and
the edges of the path $P_3$ starting in $v_3$ in each call of the recursion, and these edges are not counted multiple times
in different branches of the recursion.
Likewise, runtime of the ``merging'' part of each recursive call can be counted to the individual steps of the resulting
contraction sequence, which is of linear length.
Hence, altogether, the algorithm runs in linear time.
\end{proof}

\section{Proving the bound for bipartite planar graphs}

On a high level, the proof will follow the same path as the previous one of \Cref{lem:core}.
However, some details will be quite different, and it is thus necessary to repeat the proof in this new setting.

In a nutshell, there will be two major changes in the detailed course of the proof:
\begin{itemize}
\item Since our graph is bipartite, we will work with a plane quadrangulation (instead of a triangulation).
This has the effect that we will need to consider, in general, more than two subproblems in the recursive decomposition (actually, up to five),
and we will slightly relax the assumptions on the bounding cycle $C$ in the recursion to handle this.
However, it does not bring any big new challenges to the proof.
\item Since, again, our graph is bipartite, we immediately get that in any BFS layering, each layer is an independent set,
and so we will never create a red edge inside the same layer.
This is the crucial saving which allows us to derive a better upper bound on the red degree, and properly exploiting this saving
requires different detailed handling of the contraction steps.
\end{itemize}

Let $G$ be a bipartite graph and $r\in V(G)$ a fixed vertex.
We again consider a (fixed) BFS tree $T$ rooted at $r$, and the corresponding BFS layering $\ca L=(L_0=\{r\},L_1,L_2,\ldots)$ of~$G$.
Observe the following:
\begin{claim}[refinement of \Cref{cl:layeri3}]\label{cl:layeri3bi}
For every edge $\{v,w\}$ of bipartite $G$ with $v\in L_i$ and $w\in L_j$, we have $|i-j|=1$, and so
a contraction of a pair of vertices from $L_i$ may create new red edges only to the vertices of~$L_{i-1}\cup L_{i+1}$.
\end{claim}

The core lemma in the bipartite (quadrangulated) case now reads:
\begin{lemma}\label{lem:corebi}
Let $G$ be a plane quadrangulation and $T$ be a left-aligned BFS tree of $G$ rooted at a vertex $r\in V(G)$ of the outer face 
and defining the BFS layering $\ca L=(L_0=\{r\},L_1,L_2,\ldots)$ of~$G$.
Assume that $P_i$, $i=1,2$, are two edge-disjoint paths of $G$ of lengths at least $1$, where $P_i$ has the ends $u_i$ and $v_i$, and $u_1=u_2$.
Furthermore, assume that $P_1$ is an $r$-geodesic path and $P_2$ is a vertical path in $G$,
that $u_1=u_2\in L_\ell$ is the vertex of $P_1\cup P_2$ closest to~$r$ (i.e., at distance $\ell$ from~$r$),
and $f=\{v_1,v_2\}\in E(G)$ is an edge such that $v_1$ is to the left of $v_2$.
Then $C:=(P_1\cup P_2)+f$ is a cycle of $G$, and let $G_C$ be the subgraph of $G$ bounded by~$C$,
let $U:=V(G_C)\setminus V(P_1\cup P_2)$ and $W:=V(G)\setminus V(G_C)$.

Then there exists an $\ca L$-respecting partial contraction sequence of $G$ which contracts only pairs of vertices that are in or stem from $U$,
results in a trigraph $G^0$, and satisfies the following conditions for every trigraph $G'$ along this contraction sequence from $G$ to~$G^0$:
\begin{enumerate}[(I)']
\item For $U':=V(G')\setminus(V(P_1\cup P_2)\cup W)$ (which are the vertices that are in or stem from $U$ in~$G'$),
every vertex of $U'$ in $G'$ has red degree at most~$6$,
\label{it:max6bi}
\item the vertex $u_1$ of $P_1\cap P_2$ has no red neighbour in $U'$, 
every vertex of $P_2-u_1$ and the (two) vertices of $P_1-u_1$ at distance at most $2$ from $u_1$ have each at most $3$ red neighbours in~$U'$,
and each remaining vertex of $P_1-u_1$ has at most $5$ red neighbours~in~$U'$, 
\label{it:P12redbi}
\item at the end of the partial sequence, that is for $U^0:=V(G^0)\setminus(V(P_1\cup P_2)\cup W)$, we have \mbox{$\big|U^0\cap L_{\ell+1}[G^0]\big|\leq2$}
and $\big|U^0\cap L_j[G^0]\big|\leq1$ for all $j>\ell+1$.
\label{it:contractfinbi}
\end{enumerate}
\end{lemma}

\begin{figure}[tb]
$$
\begin{tikzpicture}[scale=0.9]
\small
\tikzstyle{every node}=[draw, shape=circle, minimum size=3pt,inner sep=0pt, fill=black]
\node at (-0.5,1) (x) {};
\node at (2,-1) (xx) {};
\node at (6,1.5) (y) {};
\node[label=above:$r$] at (3,5) (z) {};
\draw (x)--(xx)--(y) to[bend right=18] (z) (z) to[bend right=18] (x);
\draw (x)-- +(0.4,0.3); \draw (x)-- +(0.4,-0.1);
\draw (xx)-- +(-0.1,0.3); \draw (xx)-- +(0.1,0.3);
\draw (y)-- +(-0.4,0.3); \draw (y)-- +(-0.5,-0.1);
\node[label=left:$v_1$] at (1.75,0.75) (v1) {};
\node[label=right:$v_2$] at (3.75,1) (v2) {};
\node[label=right:$~{u_1=u_2}$] at (2.75,4) (u1) {};
\draw[thick, fill=gray!20] (v2) to[bend right] (u1) to[bend right] (v1) -- (v2);
\draw (z)--(3.2,4.7) node{}--(2.5,4.3) node{}--(u1);
\draw (v1)-- +(-0.2,-0.3); \draw (v1)-- +(0.2,-0.3);
\draw (v2)-- +(-0.2,-0.3); \draw (v2)-- +(0.2,-0.3);
\draw (u1)-- +(0,-0.4);
\node at (2.35,3.5){}; \node at (3.23,3.5){};
\node at (2.0,3){}; \node at (3.57,3){};
\node[label=left:$P_1~$] at (1.8,2.5){}; \node[label=right:$~P_2$] at (3.77,2.5){};
\node at (1.69,2){}; \node at (3.85,2){};
\node[label=above:$C\quad$] at (v2) {};
\node[draw=none,fill=none] at (2.75,2) {$U$};
\node[draw=none,fill=none] at (2.75,0.7) {$f$};
\node[draw=none,fill=none] at (5,1.5) {$W$};
\end{tikzpicture}
\raise18ex\hbox{\huge\boldmath~$\leadsto$~}
\begin{tikzpicture}[scale=0.9]
\small
\tikzstyle{every node}=[draw, shape=circle, minimum size=3pt,inner sep=0pt, fill=black]
\node at (-0.5,1) (x) {};
\node at (2,-1) (xx) {};
\node at (6,1.5) (y) {};
\node[label=above:$r$] at (3,5) (z) {};
\draw (x)--(xx)--(y) to[bend right=18] (z) (z) to[bend right=18] (x);
\draw (x)-- +(0.4,0.3); \draw (x)-- +(0.4,-0.1);
\draw (xx)-- +(-0.1,0.3); \draw (xx)-- +(0.1,0.3);
\draw (y)-- +(-0.4,0.3); \draw (y)-- +(-0.5,-0.1);
\node[label=left:$v_1$] at (1.75,0.75) (v1) {};
\node[label=right:$v_2$] at (3.75,1) (v2) {};
\node[label=right:$~{u_1=u_2}$] at (2.75,4) (u1) {};
\draw[thick] (v2) to[bend right] (u1) to[bend right] (v1) -- (v2);
\draw (z)--(3.2,4.7) node{}--(2.5,4.3) node{}--(u1);
\draw (v1)-- +(-0.2,-0.3); \draw (v1)-- +(0.2,-0.3);
\draw (v2)-- +(-0.2,-0.3); \draw (v2)-- +(0.2,-0.3);
\node at (2.35,3.5) (pa){}; \node at (3.23,3.5) (qa){};
\node at (2.0,3) (pb){}; \node at (3.57,3) (qb){};
\node[label=left:$P_1~$] at (1.8,2.5) (pc){}; \node[label=right:$~P_2$] at (3.77,2.5) (qc){};
\node at (1.69,2) (pd){}; \node at (3.85,2) (qd){};
%
\node[draw=none,fill=none] at (2.7,1.5) {$U^0$};
\node[draw=none,fill=none] at (5,1.5) {$W$};
\node[draw=none,fill=none,red] at (3.3,1.75) {\scriptsize$(\leq\!6)$};
\node[draw=none,fill=none,red] at (2.4,3.95) {\scriptsize$(0)$};
\node[draw=none,fill=none,red] at (1.2,2) {\scriptsize$(\leq\!5)$};
\node[draw=none,fill=none,red] at (1.5,3) {\scriptsize$(\leq\!3)$};
\node[draw=none,fill=none,red] at (4.05,3) {\scriptsize$(\leq\!3)$};
\node[draw=none,fill=none,red] at (4.3,2) {\scriptsize$(\leq\!3)$};
\node[red] at (2.7,3.5) (ra) {}; \node[red] at (2.93,3.5) (rb) {};
\node[red] at (2.8,3) (rc) {}; \node[red] at (2.6,2.5) (rd) {}; \node[red] at (2.83,2) (re) {};
\draw (u1)--(ra);
\draw[red, thick,densely dotted] (re)--(rd)--(rc)--(rb) (rc)--(ra) (qc)--(re);
\draw[red, thick,densely dotted] (pc)--(rc)--(qa) (pd)--(rd)--(qb) (ra)--(pb)--(rb) (rc)--(pa) (rd)--(pb) (re)--(pc);
\end{tikzpicture}
$$
\caption{(left) The setup of \Cref{lem:corebi}, where $v_1$ is to the left of~$v_2$.
	(right) The outcome of the claimed partial contraction sequence which contracts only vertices of $U$ inside the shaded region
	from the left and maintains bounded red degrees in the region and on its boundary $C$ formed by $P_1$ and $P_2$ and $f$.
	Not all depicted red edges really exist, and some of them may actually be black.}
\label{fig:corebilem}
\end{figure}
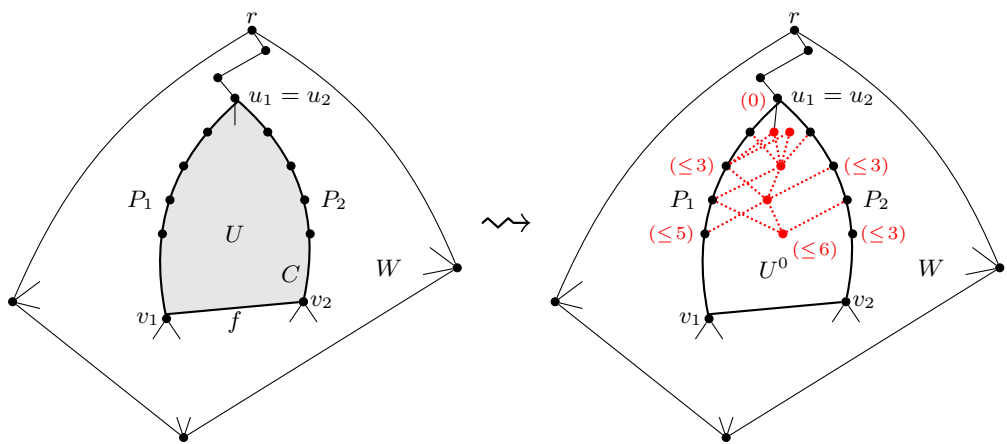

We illustrate \Cref{lem:corebi} in \Cref{fig:corebilem}.
Similarly as in the previous section, we observe that the assumptions of \Cref{lem:corebi} directly imply the following properties useful for the coming proof of \Cref{lem:corebi}.

\begin{claim}\label{clm:redalignbi}\it
Within the notation and assumptions of \Cref{lem:corebi}, for $G'$ we also have that:
\begin{enumerate}[(I)']\setcounter{enumi}{3}
\item Every red edge in $G'$ has one end in $U'$ and the other end in~$U'\cup V(P_1\cup P_2)$,
\label{it:nooutsidebi}
\item if $v\in V(P_2)\cap L_j[G']$, then there is no edge in $G'$ (red or black) from $v$ to a vertex of $(U'\cup V(P_1))\cap L_{j-1}[G']$,
\label{it:redalign2bi}
\end{enumerate}
\end{claim}
\begin{claimproof}
These properties follow in the same way as \eqref{it:nooutside} and \eqref{it:redalign2} in \Cref{clm:redalign}.
\end{claimproof}

\begin{claim}\label{clm:endreddbi}\it
With the notation and assumptions of \Cref{lem:core}, we have that if \eqref{it:contractfinbi}' is true, 
then the following holds at the end of the partial contraction sequence:
every vertex $x$ of $P_1-u_1$ has at most $2$ red neighbours in $U^0$, except when $x$ has distance $2$ from $u_1$ -- then the bound is~$3$,
and every vertex of $P_2-u_1$ has at most $1$ red neighbour in $U^0$.
\end{claim}
\begin{claimproof}
Consider a vertex $x\in V(P_1\cup P_2)\setminus\{u_1\}$ such that~$x\in L_j[G']$.
By \Cref{cl:layeri3bi}, the red neighbours of $x$ may be in $L_{j-1}[G']\cup L_{j+1}[G']$,
and by~\eqref{it:contractfinbi}' we have $|L_{j+1}[G']\cap U^0|\leq1$ and also $|L_{j-1}[G']\cap U^0|\leq1$,
except when $x$ has distance $2$ from $u_1$ and then we have $|L_{j-1}[G']\cap U^0|\leq2$.
This proves the sought upper bound ($1+1=2$ or $2+1=3$) for vertices of $P_1-u_1$.
Furthermore, if $x\in V(P_2)\setminus\{u_1\}$, then there is no such neighbour in $L_{j-1}[G']$ by \eqref{it:redalign2bi}'
and the upper bound becomes~$|L_{j+1}[G']\cap U^0|\leq1$.
\end{claimproof}

We also show how \Cref{lem:corebi} implies the first part of the main result:

\begin{proof}[Proof of \Cref{thm:twwbiplanar} (the upper bound)]
We start with a given simple bipartite planar graph $H$, and extend any plane embedding of $H$ into a simple plane quadrangulation $G$
such that $H$ is an induced subgraph of~$G$.
Then we choose a root $r$ on the outer face of $G$ and some left-aligned BFS tree of $G$ rooted at $r$ which exists by \Cref{clm:existslal}.
Let $P_1$ be the path of length $2$ starting in $r$ on the ``left'' of the outer face of $G$, and $P_2$ be the other edge incident to $r$ on the outer face of $G$.
Then $P_1$ is $r$-geodesic since $G$ is bipartite, and $P_2$ is vertical with respect to~$T$.
Now we apply \Cref{lem:corebi}.

This way we get a partial contraction sequence from $G$ to a (still bipartite by \Cref{cl:layeri3bi}) trigraph $G^0$ of maximum red degree $6$.
Similarly as in the proof of \Cref{thm:twwplanar}, we get that $L_0[G^0]=\{r\}$, $|L_1[G^0]|\leq4$, $|L_2[G^0]|\leq2$, 
and $\big|L_j[G^0]\big|\leq1$ for all $j\geq3$ by \eqref{it:contractfinbi}'.

In the final phase, we pairwise contract the vertices of $L_1[G^0]$ into one vertex, and the same with the vertices of $L_2[G^0]$, 
making $G^0$ into a path by \Cref{cl:layeri3bi}, and then naturally contract this path down to a vertex, never exceeding red degree $5<6$ along this final phase.
The restriction of the whole contraction sequence of $G$ to $V(H)$ then certifies that the twin-width of $H$ is at most~$6$, too.
\end{proof}

Now we get to the proof of core \Cref{lem:corebi} which will conclude \Cref{thm:twwbiplanar}.

\begin{proof}[Proof of \Cref{lem:corebi}]
If $U=\emptyset$, we are immediately done with the empty partial contraction sequence. So, we assume $U\not=\emptyset$.
Let $\phi$ be the bounded quadrangular face of $G_C$ incident to the edge~$f=\{v_1,v_2\}$, and $C_\phi\subseteq G$ its boundary cycle.
If $\phi$ has all vertices in $V(P_1\cup P_2)$, then let $f'$ be the edge of $C_\phi$ joining $P_1$ to $P_2$.
We then simply apply \Cref{lem:corebi} inductively to $f'$ and the corresponding subpaths of $P_1$ and~$P_2$, while the rest of the assumptions remain the same.
It immediately follows that the claimed conditions of \Cref{lem:corebi} are inherited from this inductive invocation.

Otherwise, there exist at most two vertices $v_3,v_k\in V(C_\phi)\setminus V(P_1\cup P_2)$, where $k=4$ if $|V(C_\phi)\setminus V(P_1\cup P_2)|=2$ and $k=3$ otherwise.
Up to symmetry between $v_3$ and $v_4$ if $k=4$, the counter-clockwise cyclic order of the vertices of $\phi$ is $(v_1,v_2,v_3,v_k)$.

Let $P_3$ be the (unique) vertical path of $G$ with the ends $v_3$ and $u_3$ such that $u_3\in V(P_1\cup P_2)$ and $P_3-u_3$ is disjoint from~$P_1\cup P_2$
(informally, $P_3$ is the vertical path from $v_3$ to the first vertex on~$P_1\cup P_2$).
Likewise, if $k=4$, let $P_4$ be the (unique) vertical path of $G$ with the ends $v_4$ and $u_4$ such that $u_4\in V(P_1\cup P_2\cup P_3)$ and $P_4-u_4$ is disjoint from~$P_1\cup P_2\cup P_3$.
Note that $Y:=P_1\cup P_2\cup\ldots\cup P_k$ is a tree rooted at $u_1$.

If the plane subgraph $Y\cup C_\phi$ defines only one bounded face other than $\phi$, then let $f'$ be the edge of $C_\phi$ not in $Y$.
We again simply apply \Cref{lem:corebi} inductively to $f'$ and the two subpaths of $Y$ (from the ends of $f'$ to $u_1$), which clearly satisfy the assumptions.
With the obtained partial contraction sequence, we are nearly finished and it suffices to finally contract the (one or two)
vertices of $V(C_\phi)\setminus V(P_1\cup P_2)$ -- which are not neighbours of $u_1$ in this case, within their $\ca L$-layers.
This easily preserves the claimed conditions of the lemma.

\medskip
In the remaining cases, we are going to define special subpaths of this tree, for each of the former paths.
Let $P_2^2$ be the subpath of $P_2$ from $v_2$ till the first intersection with $P_3\cup P_1$ (which may be $u_3$ or $u_1$).
Symmetrically, let $P_3^1$ be the subpath of $P_3\cup P_1$ from $v_3$ till the first intersection with $P_2$.
Let $P_3^2$ be the subpath of $P_3\cup P_2$ from $v_3$ till the first intersection with $P_4\cup P_1$ (with only $P_1$ if $k=3$)
and, symmetrically, (only if $k=4$) let $P_4^1$ be the subpath of $P_4\cup P_1$ from $v_4$ till the first intersection with $P_3\cup P_2$.
Let (again only if $k=4$) $P_4^2$ be the subpath of $P_4\cup P_3\cup P_2$ from $v_4$ till the first intersection with $P_1$.
Finally, let $P_1^1$ be the subpath of $P_1$ from $v_1$ till the first intersection with $P_2\cup\ldots\cup P_k$ (which may be $u_4$ or $u_3$ or $u_1$).
Observe that each of the paths $P_i^2$, $i\in\{2,3,4\}$, is vertical wrt.~$T$ in~$G$,
and each of the paths $P_i^1$, $i\in\{1,3,4\}$, is $r$-geodesic in~$G$.
%

For each of the (two or three) index pairs $(i,j)\in\{(2,3),(3,k),(k,1)\}$, $i\not=j$,
we consider the cycle $D_i:=(P_i^2\cup P_j^1)+\{v_i,v_j\}$.
See \Cref{fig:corebirec} (the left of), and observe that the plane region bounded by $C$ is partitioned (modulo the boundaries)
into the face $\phi$ and the regions bounded by $D_2,\ldots,D_k$.
There is one more special case we have to take care of.
If, for $i\in\{3,k\}$, we have that $u_i\in L_{\ell+1}\setminus V(P_2)$ ($u_i$ is a neighbour of $u_1$, but not on $P_2$),
and there exists the edge $\{x,y\}\in E(G)$ where $x\in L_{\ell+2}\cap V(P_i^2)$ ($x$ is the neighbour of $u_i$ on~$P_i^2$)
and $y\in L_{\ell+3}\cap V(P_j^1)$ ($y$ is the distance-$2$ neighbour of $u_i$ on $D_i$ in the other direction than $x$),
then we say that the cycle $D_i$ is {\em split} by the edge~$g_i=\{x,y\}$.
In such case, we denote by $D_i'$ the cycle of $D_i+g_i$ containing $g_i$ and not $u_i$, and by $D_i''$ the cycle of $D_i+g_i$ containing both $g_i$ and $u_i$.
Then we will consider the cycles $D_i'$ and $D_i''$ separately (instead of whole~$D_i$).

In order to unify the notation, we rename the cycles $D_2$, $D_3$ and $D_k$, after their possible splitting,
as $C_1,\ldots,C_m$ where $2\leq m\leq 2k-3$, and define the ``left'' and ``right'' subpaths $Q_1^1$, $Q_1^2$, $Q_2^1,\ldots$ of these cycles, as follows.
See on the right of \Cref{fig:corebirec}.
\begin{itemize}
\item
We let $C_1=D_2$, and then $C_2=D_3$ if $D_3$ is not split, but $C_2=D_3'$ and $C_3=D_3''$ if $D_3$ is split.
When $k=4$ (and $j\in\{3,4\}$ is the next available index),
we also let $C_j=D_4$ if $D_4$ is not split, but $C_j=D_4'$ and $C_{j+1}=D_4''$ if $D_4$ is split.
\item
For $j\in\{1,\ldots,m\}$, if $C_j=D_i$ was not split, we have $Q_j^2=P_i^2$ and $Q_j^1=P_{i'}^1$ where \mbox{$i'=(i\!\mod\!k)+1$}
(these are the same two paths as those which have defined $D_i$ above).
If $D_i$ was split by the edge $g_i$ and $C_j=D_i'$, then we have $Q_j^2:=P_i^2-u_i$ and $Q_j^1:=(P_{i'}^1\cap C_j)+g_i$.
If (again) $D_i$ was split and $C_j=D_i''$, then we take $Q_j^2$ as the single edge of $P_i^2$ incident to $u_i$,
and $Q_j^1$ as the length-$2$ subpath of $P_{i'}^1$ starting in~$u_i$.
\end{itemize}
Notice that $E(C_j)\setminus E(Q_j^1\cup Q_j^2)$ is always a single edge (which is one of the edges of $C_\phi$, or the split edge of $D_3$ or $D_4$).
Furthermore, similarly as above, we observe that each of the paths $Q_j^2$ is vertical in $G$, and each $Q_j^1$ is $r$-geodesic, for $j\in\{1,\ldots,m\}$.

\begin{figure}[t]
$$
\begin{tikzpicture}[scale=1.2]
\small
\tikzstyle{every node}=[draw, shape=circle, minimum size=3pt,inner sep=0pt, fill=black]
\node[label=left:$v_1$] at (1,0.5) (v1) {};
\node[label=right:$v_2$] at (5,1) (v2) {};
\node[label=right:$~~{u_1=u_2}$] at (2.75,5) (u1) {};
\node[label=below:$v_3$] at (3.75,1.4) (v3) {};
\node[label=below:$v_4$] at (2,1.2) (v4) {};
\draw[draw=none, fill=gray!20] (v1)--(v4)--(2.75,5) to[bend right=44] (v1);
\draw[draw=none, fill=gray!20] (v4)--(v3)--(2.75,5) to[bend right=44] (v4);
\draw[draw=none, fill=gray!20] (v3)--(v2) to[bend right=44] (2.75,5)--(v3);
\draw[thick] (v2) to[bend right=44] (2.75,5) to[bend right=44] (v1) (v1) -- (v2);
\draw (v2)--(v3)--(v4)--(v1);
\node[label=right:$~u_3$] at (4.35,4) (u3) {};
\node[label=left:$u_4~$] at (2.1,4.63) (u4) {};
\draw (v3) to[bend left] (u3) (v4) to[bend right=11] (u4);
\draw (u1)-- +(-0.2,0.3) (u1)-- +(0.1,-0.4);
\draw (v1)-- +(-0.2,-0.3); \draw (v1)-- +(0.2,-0.3);
\draw (v2)-- +(-0.2,-0.3); \draw (v2)-- +(0.2,-0.3);
\node[draw=none,fill=none, label=left:$P_1~$] at (0.85,3.5){};
\node[draw=none,fill=none, label=right:$~P_2$] at (5,3.5){};
\node[draw=none,fill=none] at (3.85,2) {$P_3^1$};
\node[draw=none,fill=none] at (4.85,2) {$P_2^2$};
\node[draw=none,fill=none] at (3.9,3.8) {$P_3^2$};
\node[draw=none,fill=none] at (2.4,4.2) {$P_4^1$};
\node[draw=none,fill=none] at (1.95,2) {$P_4^2$};
\node[draw=none,fill=none] at (0.95,2) {$P_1^1$};
\node[draw=none,fill=none] at (4.3,3) {$D_2$};
\node[draw=none,fill=none] at (1.6,3.3) {$D_4$};
\node[draw=none,fill=none] at (3,3.4) {$D_3$};
\node[draw=none,fill=none] at (2.8,1) {$\phi$};
\end{tikzpicture}
\qquad\qquad
%
\begin{tikzpicture}[scale=1.2]
\small
\tikzstyle{every node}=[draw, shape=circle, minimum size=3pt,inner sep=0pt, fill=black]
\node[label=left:$v_1$] at (1,0.5) (v1) {};
\node[label=right:$v_2$] at (5,1) (v2) {};
\node[label=right:$~~{u_1=u_2}$] at (2.75,5) (u1) {};
\node[label=below:$v_3$] at (3.75,1.4) (v3) {};
\node[label=below:$v_4$] at (2,1.2) (v4) {};
\draw[draw=none, fill=gray!20] (v1)--(v4)--(2.75,5) to[bend right=44] (v1);
\draw[draw=none, fill=gray!20] (v4)--(v3)--(2.75,5) to[bend right=44] (v4);
\draw[draw=none, fill=gray!20] (v3)--(v2) to[bend right=44] (2.75,5)--(v3);
\draw[thick] (v2) to[bend right=44] (2.75,5) to[bend right=44] (v1) (v1) -- (v2);
\draw (v2)--(v3)--(v4)--(v1);
\node[label=right:$~u_3$] at (4.35,4) (u3) {};
\node[label=left:$u_4~$] at (2.1,4.63) (u4) {};
\node[label=right:$~x$] at (2.2,4) (xx) {}; \node at (1.45,4) (xxx) {};
\node[label=left:$y~$] at (1.12,3.5) (yy) {};
\draw (v3) to[bend left] (u3) (v4) to[bend right=11] (u4) (xx) to[bend left=3] (yy);
\draw (u1)-- +(-0.2,0.3) (u1)-- +(0.1,-0.4);
\draw (v1)-- +(-0.2,-0.3); \draw (v1)-- +(0.2,-0.3);
\draw (v2)-- +(-0.2,-0.3); \draw (v2)-- +(0.2,-0.3);
\node[draw=none,fill=none, label=left:$P_1~~$] at (0.85,3.5){};
\node[draw=none,fill=none, label=right:$~P_2$] at (5,3.5){};
\node[draw=none,fill=none] at (3.85,2) {$Q_1^1$};
\node[draw=none,fill=none] at (4.85,2) {$Q_1^2$};
\node[draw=none,fill=none] at (3.9,3.8) {$Q_2^2$};
\node[draw=none,fill=none] at (2.5,3.3) {$Q_2^1$};
\node[draw=none,fill=none] at (2.02,2.8) {$Q_3^2$};
\node[draw=none,fill=none] at (1.2,3.2) {$Q_3^1$};
\node[draw=none,fill=none] at (2.4,4.4) {$Q_4^2$};
\node[draw=none,fill=none] at (1.2,4.2) {$Q_4^1$};
\node[draw=none,fill=none] at (4.3,2.8) {$C_1\!=\!D_2$};
\node[draw=none,fill=none] at (1.4,2.2) {$C_3\!=\!D'_4$};
\node[draw=none,fill=none] at (3,2.4) {$C_2\!=\!D_3$};
\node[draw=none,fill=none] at (1.8,4.04) {$D''_4$};
\node[draw=none,fill=none] at (1.8,3.6) {$g_4$};
\end{tikzpicture}
$$
\caption{(left) An illustration of the general setup of the main induction step in the proof of \Cref{lem:corebi};
	the region bounded by $C$ is ``divided'' into the face $\phi$ and up to three subregions bounded by the cycles $D_2$, $D_3$ and $D_4$
	(which are in turn defined by the vertical paths $P_3$ and~$P_4$ inside~$C$).
	(right) The division from the left is further refined at the cycle $D_4$ which is split by the edge $g_4$,
	the defining cycles of this refined division are further renamed as $C_1=D_2$, $C_2=D_3$, $C_3=D_4'$ and $C_4=D_4''$,
	and the paths $Q_j^1$ and $Q_j^2$ of each $C_j$, $j\in\{1,2,3,4\}$, are outlined in the picture.}
\label{fig:corebirec}
\end{figure}
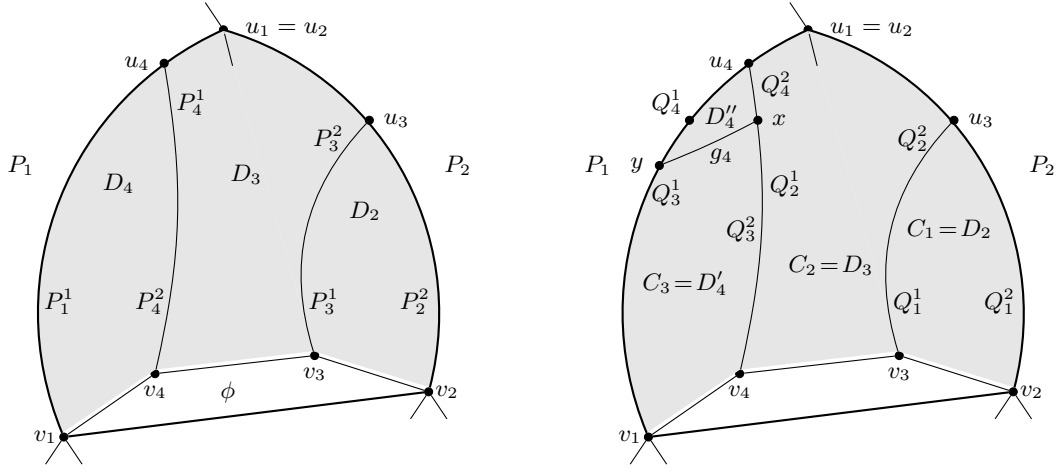

\medskip
The key conclusion from the initial step of this proof is that we can apply \Cref{lem:corebi} inductively (and independently)
to each of the cycles $C_j$, $j\in\{1,\ldots,m\}$, more specifically to the pair of paths $Q_j^1$, $Q_j^2$.
However, the application is not as straightforward as in the proof of \Cref{lem:core}, since we have more than two subcases
in general and we are going to merge them pairwise; the outcome of the subcase of $C_1$ with that of $C_2$, then the result
with the outcome of the subcase $C_3$, then with~$C_4$\,\dots
To formally manage this process, via nested induction, we are going to show that the merged region
(considering it with a generalized ``bounding'' which is to be defined below) will again satisfy the conditions 
\eqref{it:max6bi}', \eqref{it:P12redbi}' and \eqref{it:contractfinbi}' of \Cref{lem:corebi}, and consequently also
the additional conditions of \Cref{clm:redalignbi}.

We say that a subgraph $H_1\subseteq G$ is {\em properly bounded} by paths $R^1,R^2\subseteq H_1$
(with respect to implicit~$G$, and in this order of $R^1$, $R^2$), if the following holds.
The path $R^i$, $i=1,2$, has the ends $s_i$ and $t_i$ and $s_1=s_2$ is the vertex of $H_1$ closest to~$r$ in~$G$,
and $R^1$ is an $r$-geodesic path and $R^2$ is a vertical path in $G$.
Moreover, there is no edge between a vertex of $V(H_1)\setminus V(R^1\cup R^2)$ and a vertex of $V(G)\setminus V(H_1)$.
Notice that, under the assumptions of \Cref{lem:corebi}, the subgraph $G_C$ is properly bounded by $P_1$ and~$P_2$,
however, here we do not explicitly require the existence of an edge~$\{t_1,t_2\}$.
We shortly say that the vertex $s_1=s_2$ is the {\em sink} of $H_1$, 
and that the vertices of $V(H_1)\setminus V(R^1\cup R^2)$ form the {\em interior} of $H_1$, 
those of $V(G)\setminus V(H_1)$ the {\em exterior} of $H_1$ and those of $V(R^1\cup R^2)$ the {\em boundary} of~$H_1$.

The rest of the proof relies on the following technical claim:

\begin{claim}\label{clm:subcorebi}
Assume a subgraph $H_0\subseteq G$ that is properly bounded by paths $R_0^1,R_0^2\subseteq H_0$,
together with two subgraphs $H_i\subseteq H_0$, $i=1,2$, properly bounded by paths $R_i^1,R_i^2$, 
and an $r$-geodesic path $R'$ in $G$, such that the~following~hold:
\begin{enumerate}[a)]
\item $V(R')=V(H_1)\cap V(H_2)$, and $V(R')\setminus V(R_0^1\cup R_0^2)\subseteq V(R_1^2)\cap V(R_2^1)$,
\label{it:Rmid}
\item $V(R_0^1)\subseteq V(R_1^1)\cup V(R_2^1)$ and $V(R_0^2)\subseteq V(R_1^2)\cup V(R_2^2)$,
\label{it:Rsides}
\item there is at most one vertex in $V(R')\cap V(R_0^2)$ and then it is the sink of $H_2$,
and for $z$ being the end of $R_1^1$ farther from $r$, there is at most one vertex in $V(R')\cap V(R_0^1-z)$ and then it is the sink of $H_1$,
\label{it:Sinks}
\item if $z\in V(R')\cap V(R_0^1)$, then the length of $R_1^1$ is exactly $2$ and $z$ is at distance $3$ from the sink of~$H_0$ (cf.~the split case in the above definitions).
\label{it:Split}
\end{enumerate}
\smallskip
If, for $i=1,2$, there is an $\ca L$-respecting partial contraction sequence $\tau_i$ of $G$ contracting only among the interior vertices of $H_i$,
such that $\tau_i$ satisfies the conditions \eqref{it:max6bi}', \eqref{it:P12redbi}' and \eqref{it:contractfinbi}' of \Cref{lem:corebi}
(specifically, with $P_1$ and $P_2$ replaced by $R_i^1$ and $R_i^2$, respectively, and $U$ and $W$ being the interior and the exterior vertices of~$H_i$),
then there is an $\ca L$-respecting partial contraction sequence $\tau_0$ of $G$ contracting only among the interior vertices of $H_0$,
such that $\tau_0$ satisfies the conditions \eqref{it:max6bi}', \eqref{it:P12redbi}' and \eqref{it:contractfinbi}' of \Cref{lem:corebi}, too.
\end{claim}

We first show how \Cref{clm:subcorebi} finishes the proof of \Cref{lem:corebi}.
Let $G_j\subseteq G$, $j\in\{1,\ldots,m\}$, denote the subgraph bounded by the cycle $C_j$ in $G$.
We apply \Cref{clm:subcorebi} to $H_0=G_2':=G_1\cup G_2$ and $H_1=G_1$, $H_2=G_2$, where the partial contraction sequences
$\tau_i$, $i\in\{1,2\}$, assumed by the claim, are obtained for $i=1,2$ by an inductive invocation of \Cref{lem:corebi} to $C_i$
with the above defined paths $Q_i^1=R_i^1$, $Q_i^2=R_i^2$.
The remaining conditions of \Cref{clm:subcorebi} can be routinely verified from the definition of the cycles $C_1$ and~$C_2$.
In particular, $H_0=G_2'$ is properly bounded by appropriate paths $R_0^1,R_0^2$ since any edge between an interior and
an exterior vertex of $H_0$ would have to cross the face $\phi$ of the plane graph $G$ which is absurd.

Subsequently, for $j=3,\ldots,m$, we analogously apply \Cref{clm:subcorebi} to $H_0=G_j':=G_{j-1}'\cup G_j$ and $H_1=G_j$, $H_2=G_{j-1}'$,
where the partial contraction sequence $\tau_2$ has been obtained inductively in the previous iteration for $G_{j-1}'$,
and $\tau_1$ is obtained by an inductive invocation of \Cref{lem:corebi} to $C_j$ with the above defined paths $Q_j^1=R_1^1$, $Q_j^2=R_1^2$.
Since $G_m'=G_C$, as considered by \Cref{lem:corebi}, the proof of \Cref{lem:corebi} is finished by this.

\medskip
Now it remains to prove \Cref{clm:subcorebi}.
We first clarify that both \Cref{clm:redalignbi} and \Cref{clm:endreddbi} remain true in the extended setting of \Cref{clm:subcorebi}.
For property \eqref{it:nooutsidebi}' of \Cref{clm:redalignbi}, this follows from the definition of $H_0$ being properly bounded,
and for property \eqref{it:redalign2bi}' the arguments of \Cref{clm:redalignbi} stay the same since $R_0^2$ is vertical in~$G$.
The arguments of \Cref{clm:endreddbi} also stay the same.


We apply the $\ca L$-respecting partial contraction sequences $\tau_2$ and $\tau_1$ on $G$ in this order (i.e., $\tau_1$ after whole~$\tau_2$).
This is sound since the subgraphs $H_1$ and $H_2$ are properly separated, and so no contraction in $\tau_2$ has an effect
in the neighbourhood of the interior vertices of $H_1$ and vice versa.
Let $\tau'=\tau_2.\tau_1$ be the concatenated sequence contracting from original $G$ to the trigraph $G^2$.
We first verify that the conditions \eqref{it:max6bi}' and \eqref{it:P12redbi}' are satisfied along $\tau'$ (while \eqref{it:contractfinbi}' is void for $G^2$ at this stage).
This is immediate for \eqref{it:max6bi}' and the interior vertices of $H_1$ and of $H_2$ since they do not get any other red neighbour than from their partial sequence.
Each vertex $x\in V(R')\setminus V(R_0^1\cup R_0^2)$ satisfies \eqref{it:max6bi}' along $\tau_2$,
using \eqref{it:P12redbi}' for $H_2$ by \Cref{clm:subcorebi}.\ref{it:Rmid},
and then (still in $\tau_2$) it ends up with at most $1$ red neighbour in the (contracted) interior of $H_2$ by \Cref{clm:endreddbi} for $H_2$.
Subsequently, $x$ never exceeds red degree of $5+1=6$ along $\tau_1$ thanks to again \Cref{clm:subcorebi}.\ref{it:Rmid}\, and \eqref{it:P12redbi}' for $H_1$.

We are left with the vertices $y\in V(R_0^1\cup R_0^2)$ and the condition \eqref{it:P12redbi}' for $H_0$.
Assume first $y\not\in V(H_1)\cap V(H_2)$.
By \Cref{clm:subcorebi}.\ref{it:Rsides}, validity of the condition \eqref{it:P12redbi}' is inherited from its validity for $H_1$ or $H_2$;
this is immediate except possibly for the case of $y\in V(R_0^1)$ at distance at most $2$ from the sink of $H_0$,
in which it suffices to observe that the sinks of $H_1$ and $H_2$ cannot be farther from $y$ than the sink of~$H_0$.

Now, let $y\in V(H_1)\cap V(H_2)=V(R')$. In this case, $y$ is the sink of $H_1$ or $H_2$ by \Cref{clm:subcorebi}.\ref{it:Sinks}, 
and then it receives no red edge from the other sequence and \eqref{it:P12redbi}' holds.
Alternatively, we may have the special case of $y=z$ (the far end of~$R'$) and $y\in V(R_0^1)$, 
in which the assumption of \Cref{clm:subcorebi}.\ref{it:Split} applies.
Moreover, in this case the sink of $H_2$ is identical to the sink of $H_0$, and so at the end of $\tau_2$
the vertex $y$ has at most $2$ red neighbours in the (contracted) interior of $H_2$ by \Cref{clm:endreddbi} for~$H_2$.
Subsequently, along $\tau_1$, the number of red neighbours of $y$ in the interior of $H_1$ always stays at most~$3$ by
\eqref{it:P12redbi}' for $H_1$ (since $y$ is at distance $2$ from the sink of $H_1$).
Altogether, this means at most $2+3=5$ red neighbours in the interior of $H_0$ along whole $\tau'$, as needed in \eqref{it:P12redbi}' for $H_0$.

Verification of the conditions \eqref{it:max6bi}' and \eqref{it:P12redbi}' along $\tau'$ is finished.

\medskip
We continue contractions from the trigraph $G^2$ to $G^3$ along a special $\ca L$-respecting partial contraction sequence $\tau_3$ as follows.
Let $W$ denote the exterior vertices of $H_0$, and $U':=V(G^2)\setminus(V(R_0^1\cup R_0^2)\cup W)$ be the vertices of $G^2$
that stem by contractions from the interior of~$H_0$.
If the sink $s_0$ of $H_0$ is identical to the sinks of $H_1$ and of $H_2$, and $s_0\in L_\ell$, then we have up to $5$ vertices in
$U'\cap L_{\ell+1}[G^2]$, and we iteratively (in any order) contract those pairs of them with the same adjacency relation to $s_0$.
This preserves validity of the condition \eqref{it:P12redbi}' for $H_0$, and regarding \eqref{it:max6bi}' we use \Cref{cl:layeri3bi},
according to which the red neighbours of the contracted vertices from $U'\cap L_{\ell+1}[G^2]$ are only in the layer $\ell+2$
and we clearly have $|(V(G^2)\setminus W)\cap L_{\ell+2}[G^2]|\leq5$ by \eqref{it:contractfinbi}' for $H_1$ and~$H_2$.
If the sink $s_2$ of $H_2$ belongs to $V(R_0^2-s_0)$, and $s_2\in L_p$, then we contract the up to two vertices of $U'\cap L_{p+1}[G^2]$ into one.
This adds one red edge to $s_2$, but this is safe by \Cref{clm:endreddbi} for $H_1$ since $1+1<3$.
Again, validity of the conditions \eqref{it:max6bi}' and \eqref{it:P12redbi}' along $\tau_3$ is easily preserved.
In all other cases we leave $\tau_3$ empty.

Then we construct an $\ca L$-respecting partial contraction sequence $\tau_4$ from $G^3$ to $G^4$.
Let $U_0:=V(G^3)\setminus(V(R_0^1\cup R_0^2)\cup W)$ be the vertices of $G^3$ that are or stem by contractions from the interior of~$H_0$.
Let $U_3:=U_0\cap V(R')$, and for $i=1,2$ let $U_i\subseteq U_0$ be the vertices that are or stem by contractions from the interior of~$H_i$.
For an informal illustration, $(U_1,U_2,U_3)$ is ``nearly'' a partition of~$U_0$, except that $U_0\cap L_{\ell+1}[G^3]$ may belong to both $U_1$ and $U_2$ due to the contractions in~$\tau_3$.
Note that $\big|U_2\cap L_{p}[G^3]\big|\leq1$ for all $p\geq\ell+2$, by $\tau_3$, and similarly $\big|U_1\cap L_{p}[G^3]\big|\leq1$
for all $p\geq\ell+2$ except when $p$ is the layer next to the sink of $H_1$ (then it is $\cdot\leq2$),
and $\big|U_3\cap L_{p}[G^3]\big|\leq1$ since $R'$ is $r$-geodesic.
Furthermore, $U_3\cap L_{\ell+1}[G^3]=\emptyset$ as ensured by the sequence~$\tau_3$.

The sequence $\tau_4$ is going to contract each layer of $U_2\cup U_3$ into one vertex, starting from the farthest layer till layer~$\ell+2$.
Formally, let $d:=\max\{i: (U_2\cup U_3)\cap L_{i}[G^3]\not=\emptyset\}$.
For $i=d,d-1,\ldots,\ell+2$ in this order, if $\big|(U_2\cup U_3)\cap L_{i}[G^3]\big|>1$, then $\tau_4$ contracts these two vertices in layer $i$.
That is all for $\tau_4$, and we now verify the conditions \eqref{it:max6bi}' and \eqref{it:P12redbi}' for any trigraph $G'$ along $\tau_4$.
Let $H_0':=G'-W$ be the corresponding trigraph contracted from $H_0$, 
and let~$U_3'$ be the vertices of $G'$ that stem by $\tau_4$-contactions from~$U_2\cup U_3$.
We consider all possibly affected vertices of $H_0'$ as follows:
\begin{itemize}
\item The sink $s_0$ of $H_0$ (which is the sink of $H_0'$, too) has got no red edge.
\item For $x\in V(R_0^2)\setminus\{s_0\}$, where $x\in L_{p}[G']$, we have no red neighbours in $L_{p-1}[G']$ by \eqref{it:redalign2bi}' applied to~$H_2$.
In $L_{p+1}[G']$, there is always at most one red neighbour of $x$ by the contraction process of $\tau_3.\tau_4$, 
except when $x$ is the sink of $H_2$ and a second red neighbour of $x$ may exist in $U_1$.
Hence, \eqref{it:P12redbi}' is true for~$x$.

\item For $x\in V(R_0^1)\setminus\{s_0\}$, where $x\in L_{p}[G']$, we distinguish two subcases.
If $x\in V(H_2)$, then $\tau_4$ does not contract in the neighbourhood of $x$ (and there is nothing to solve), except
when $x$ is also the sink of $H_1$ -- then $x$ has at most two (two if $p=\ell+2$) red neighbours in $U_2\cap L_{p-1}[G']$ 
and an additional one in contracted $U_3'\cap L_{p+1}[G']$.

Secondly, if $x\in V(H_1)\setminus V(H_2)$, then $x$ has at most two red neighbours in layer $p-1$,
that is one in $U_1\cap L_{p-1}[G']$ and one possibly in $U_3'\cap L_{p-1}[G']$.
In addition to that, $x$ may have one red neighbour in $U_1\cap L_{p+1}[G']$, but no red neighbour in $(U_2\cup U_3)\cap L_{p+1}[G']$
and no such in $U_3'\cap L_{p+1}[G']$ -- the latter is due to \eqref{it:redalign2bi}' applied to~$H_1$.

In both subcases, we have that \eqref{it:P12redbi}' is true for~$x$.

\item If $x\in U_1\cap L_{p}[G']$ and $p\geq\ell+2$, then $x$ cannot be adjacent to vertices in $U_2\cup V(R_0^2)$ since $H_1$ was properly separated.
There are at most $4$ red neighbours of $x$ in $(V(R_0^1)\cup U_1\cup U_3\cup U_3')\cap L_{p-1}[G']$ (two of which may be in $U_1\cap L_{p-1}[G']$).
Additionally, there are up to $2$ red neighbours of $x$ in $(V(R_0^1)\cup U_1)\cap L_{p+1}[G']$, 
but $x$ cannot be a neighbour of the vertex of $U_3'\cap L_{p+1}[G']$ due to \eqref{it:redalign2bi}' applied to~$H_1$.
We have that \eqref{it:max6bi}' is true for~$x$.

\item If $x\in(U_1\cup U_2)\cap L_{\ell+1}[G']$, then $x$ is not contracted by $\tau_4$, and $x$ has potential red neighbours
only in $L_{\ell+2}[G']\cap V(H_0')$ of at most $6$ vertices.
So, \eqref{it:max6bi}' is true for~$x$.

\item It remains to consider $x\in(U_2\cup U_3\cup U_3')\cap L_{p}[G']$ where $p\geq\ell+2$.
Since $H_2$ is properly separated, and by \eqref{it:redalign2bi}' applied to~$H_1$, there is no red edge from $x$
to $(V(R_0^1)\cup U_1)\cap L_{p-1}[G']$.
For the rest, there are hence at most $4$ red neighbours of $x$ in $(U_2\cup U_3\cup U_3')\cap L_{p-1}[G']$.
There is also no red edge from $x$ to $V(R_0^2)\cap L_{p+1}[G']$ by \eqref{it:redalign2bi}' applied to~$H_2$.
If $x\not\in U_3'$, there are at most $2$ red neighbours of $x$ in layer $p+1$, and otherwise there can be up to
$3$ red neighbours of $x$ there in $(V(R_0^1)\cup U_1\cup U_3')\cap L_{p+1}[G']$.

The latter subcase indicates a possible problem of $x\in U_3'$ reaching the red degree of $4+3>6$, but for this to actually happen,
it would have to be $|U_2\cap L_{p-1}[G']|=2$, meaning that $p=\ell+2$, and $x$ would need to have a red edge to
$V(R_0^1)\cap L_{p+1}[G']$, meaning that there was such a black edge in $G^3$ and so already in~$G$.
This would be the case of a split cycle $D$ resolved in the initial step of the proof by another decomposition,
and hence this cannot be true here.
Again, \eqref{it:max6bi}' holds for~$x$.
\end{itemize}

We continue from the trigraph $G^4$ to $G^5$ with an $\ca L$-respecting partial contraction sequence $\tau_5$ as follows.
We recall the set $U_1$ from the previous phase which stays the same in $G^4$, 
and the set $U_3'$ as those vertices of $G^4$ which are in or stem by contractions from $U_2\cup U_3$ of~$G^3$.
The purpose of $\tau_5$ is similar to that of the previous special sequence $\tau_3$; to contract the possible pair of vertices
of $U_1$ which are neighbours of the sink of $H_1$.
Let $s_1$ denote the sink of $H_1$, and assume $s_1\in L_{p}[G^4]$.
If $p=\ell$ or $|U_1\cap L_{p+1}[G^4]|\leq1$, then $\tau_5$ is empty.
Otherwise, if $p=\ell+1$, then $\tau_5$ simply contracts the vertex pair in $U_1\cap L_{p+1}[G^4]$.
If $p>\ell+1$, then $\tau_5$ contracts the single vertex $z\in U_3'\cap L_{p+1}[G^4]$ 
with one and then with the other one of the two vertices of $U_1\cap L_{p+1}[G^4]$.

We verify validity of the conditions \eqref{it:max6bi}' and \eqref{it:P12redbi}' along $\tau_5$.
In the case of $p=\ell+1$, the only vertex with possibly increased red degree is $s_1$.
However, in this case $s_1$ has no red neighbour in layer $p-1=\ell$ and at most one in layer $p+1$ by the analysis of $\tau_4$.
After $\tau_5$, the red degree of $s_1$ is thus at most~$1+1=2$ which is good for \eqref{it:P12redbi}'.
In the case of $p>\ell+1$, the potential red neighbours of the vertices being contracted with $z$ 
are the at most three vertices of $V(R_0^1\cup R_0^2\cup U_3')\cap L_{p-1}[G^4]$ plus 
the at most three vertices of $V(R_0^1\cup U_1\cup U_3')\cap L_{p+1}[G^4]$ (recall the previous analysis of $\tau_4$), 
which sums to the red degree at most~$6$ along $\tau_5$ and is good for \eqref{it:max6bi}'.
Furthermore, for each of these considered six vertices, $z$ has been counted as their red neighbour in
the previous analysis of $\tau_4$, and so their maximum red degree does not increase along $\tau_5$.

\medskip
At last, we construct an $\ca L$-respecting partial contraction sequence $\tau_6$ from $G^5$ to $G^6$.
We denote by $U_4$ the set of vertices of $G^5$ including $U_1\cap V(G^5)$ and the possible vertex resulting by the contractions of $\tau_5$,
and by $U_5:=V(G^5)\setminus(V(R_0^1\cup R_0^2)\cup W\cup U_4)$ the set of those vertices of $G^5$ which 
stem by contractions from $U_2\cup U_3$ of~$G^3$ in $\tau_4$ and are not touched by~$\tau_5$.
The sequence $\tau_5$ is going to contract each layer of $U_4\cup U_5$ into one vertex, starting from the farthest layer till layer~$\ell+2$.
Formally, let $d:=\max\{i: (U_4\cup U_5)\cap L_{i}[G^5]\not=\emptyset\}$.
For $i=d,d-1,\ldots,\ell+2$ in this order, if $\big|(U_4\cup U_5)\cap L_{i}[G^5]\big|>1$,
then $\tau_5$ contracts the pair of vertices of $(U_4\cup U_5)\cap L_{i}[G^5]$.

We now verify the conditions \eqref{it:max6bi}' and \eqref{it:P12redbi}' for any trigraph $G'$ along $\tau_6$.
Let~$U_5'$ denote the set of vertices of $G'$ that stem by $\tau_6$-contactions from~$U_4\cup U_5$.
We perform an analogous case analysis of all possibly affected vertices:
\begin{itemize}
\item The sink $s_0$ of $H_0$ has again got no red edge.
\item For $x\in V(R_0^2)\setminus\{s_0\}$, where $x\in L_{p}[G']$, we have no red neighbours in $L_{p-1}[G']$ by \eqref{it:redalign2bi}',
and there are at most $2$ candidates for a red neighbour of $x$ in $U_5'\cap L_{p+1}[G']$ (this is after a $\tau_6$-contaction in layer~$p+1$).
Hence, \eqref{it:P12redbi}' is true for~$x$.

\item For $x\in V(R_0^1)\setminus\{s_0\}$, where $x\in L_{p}[G']$, and after a $\tau_6$-contaction in layer~$p+1$,
the vertex $x$ has at most one red neighbour in layer $p+1$, in $U_5'\cap L_{p+1}[G']$.
In layer $p-1$, there is possibly one red neighbour in $(U_5\cup U_5')\cap L_{p-1}[G']$ and one in $U_4\cap L_{p-1}[G']$.
Together, \eqref{it:P12redbi}' is true for~$x$.

\item If $x\in(U_1\cup U_3')\cap L_{\ell+1}[G']$, then $x$ is not contracted by $\tau_5$, and \eqref{it:max6bi}' is true for~$x$ as previously.

\item It remains to consider $x\in(U_4\cup U_5\cup U_5-')\cap L_{p}[G']$ where $p\geq\ell+2$, and after a $\tau_6$-contaction in layer~$p+1$.
The only possible red neighbours of $x$ in layer $p+1$ are in $(V(R_0^1)\cup U_5')\cap L_{p+1}[G']$ which has at most $2$ vertices.
There is again no red edge from $x$ to $V(R_0^2)\cap L_{p+1}[G']$ by property \eqref{it:redalign2bi}'.
In layer $p-1$, there are altogether at most $4$ vertices in $(V(R_0^1\cup R_0^2)\cup U_5\cup U_5')\cap L_{p-1}[G']$.
So, the bound on red degree at most $6$ holds and \eqref{it:P12redbi}' is true for~$x$.
\end{itemize}

The last step in the proof of \Cref{clm:subcorebi} is to observe that the concatenated sequence 
$\tau_0:=(\tau'.\tau_3.\tau_4.\tau_5.\tau_6)$, with the outcome $G^6$ as above,
satisfies also the condition \eqref{it:contractfinbi}'.
\Cref{lem:corebi} is thus proved.
\end{proof}

\begin{proof}[Proof of \Cref{thm:twwbiplanar} (the algorithmic part)]
We can construct a simple plane quadrangulation $G\supseteq H$ in linear time using standard planarity algorithms,
and then construct a left-aligned BFS tree $T\subseteq G$ again in linear time by \Cref{clm:existslal}.
In the rest, we straightforwardly implement the recursive division of $G$ as used in the proof of \Cref{lem:corebi}
and construct the contraction sequence of $H$ on return from the recursive calls as defined in the proof
of the lemma, specifically in \Cref{clm:subcorebi}.
Note that we again do not need at all to construct the intermediate trigraphs along the constructed contraction sequence.
Hence this proof is the same as that of \Cref{thm:twwplanar}.
\end{proof}

\section{Conclusion}

We have significantly improved the best upper bound on the twin-width of planar graphs,
and we seem to be approaching its true maximum value on planar graphs.
While some arguments in the proof of \Cref{lem:core} are still loose and can be further improved by one or two down,
there are several specific situations in the proof in which one cannot go lower than $9$, e.g., in the ``game ending'' trigraphs in \Cref{fig:endcases}.
Though, we have no specific conjecture about what the true maximum can be.

We have also given an upper bound on the twin-width of bipartite planar graphs, which have not yet been
specifically targetted by twin-width research.
The bound of $6$ again seems to approach the true maximum value on bipartite planar graphs and,
in particular, there does not seem to be room for further improvement in the proof of \Cref{lem:corebi}.

The natural next step of our research is to have a closer look at suitable lower-bound planar examples.
Another natural direction of research is to extend this result to other graph classes which are related
to planar graphs, and especially to those classes for which explicit twin-width upper bounds have already
been given using the product-structure machinery (e.g., for $k$-planar graphs).

\bibliography{tww}

\end{document}